\documentclass{article}
\pdfoutput=1 
\usepackage[utf8]{inputenc}
\usepackage[british]{babel}
\usepackage[T1]{fontenc}

\usepackage[all]{xy}
\usepackage[a4paper,headheight=3ex,body={15cm,23cm}]{geometry}
\usepackage{lmodern, graphicx, xcolor, mathtools, amssymb, amsthm, amsfonts, booktabs, enumitem, hyperref}

\hypersetup{
    colorlinks=true,
    pdftitle={Mordell-Weil groups as Galois modules},
    pdfauthor={Christian Wuthrich},
    breaklinks=true
}

\usepackage[textsize=tiny, textwidth=20mm]{todonotes}

\definecolor{dblue}{rgb}{0,0,0.7}
\newtheoremstyle{mythm}{2ex}{2ex}{\color{dblue}}{}{\bf\color{dblue}}{.}{ }{}
\theoremstyle{mythm}
\newtheorem{thm}{Theorem}
\newtheorem{prop}[thm]{Proposition}
\newtheorem{cor}[thm]{Corollary}
\newtheorem{lem}[thm]{Lemma}

\theoremstyle{definition}
\newtheorem{assu}{Assumption}

\definecolor{dgreen}{rgb}{0,0.5,0}
\newcommand{\hl}[1]{\textcolor{dgreen}{#1}}

\renewcommand{\geq}{\geqslant}
\renewcommand{\leq}{\leqslant}

\DeclareMathOperator{\No}{N}
\DeclareMathOperator{\coker}{coker}
\DeclareMathOperator{\Gal}{Gal}
\DeclareMathOperator{\Ext}{Ext}
\DeclareMathOperator{\Hom}{Hom}
\DeclareMathOperator{\Ind}{Ind}
\DeclareMathOperator{\res}{Res}
\DeclareMathOperator{\rk}{rk}
\DeclareMathOperator{\ord}{ord}
\DeclareMathOperator{\Aut}{Aut}
\DeclareMathOperator{\GL}{GL}
\DeclareMathOperator{\SL}{SL}
\DeclareMathOperator{\sign}{sign}

\newcommand{\Sel}{\mathcal{S}}

\newcommand{\Selc}{\mathfrak{S}}
\newcommand{\Relc}{\mathfrak{R}}
\newcommand{\Ens}{\tilde{\mathcal{E}}}

\newcommand{\FF}{\mathbb{F}}
\newcommand{\QQ}{\mathbb{Q}}
\newcommand{\ZZ}{\mathbb{Z}}

\newcommand{\CC}{\mathbb{C}}

\newcommand{\QZp}{{}^{\QQ_p}\!/\!{}_{\ZZ_p}}
\newcommand{\cyclic}[1]{{}^{\ZZ}\!/\!{}_{#1\ZZ}}
\newcommand{\hatox}{\mathbin{\hat\otimes}}
\newcommand{\twist}[1]{\breve{#1}}

\DeclareFontFamily{U}{wncy}{}
\DeclareFontShape{U}{wncy}{m}{n}{<->wncyr10}{}
\DeclareSymbolFont{mcy}{U}{wncy}{m}{n}
\DeclareMathSymbol{\Sha}{\mathord}{mcy}{"58}

\begin{document}
  \title{Mordell-Weil groups as Galois modules}
  \author{Thomas Vavasour and Christian Wuthrich}
  \maketitle

  \begin{abstract}
    We study the action of the Galois group $G$ of a finite extension $K/k$ of number fields on the points on an elliptic curve $E$.
    For an odd prime $p$, we aim to determine the structure of the $p$-adic completion of the Mordell-Weil group $E(K)$ as a $\ZZ_p[G]$-module only using information of $E$ over $k$ and the completions of $K$.
  \end{abstract}

\setcounter{tocdepth}{1}
\tableofcontents

\section{Introduction}

Let $E$ be an elliptic curve defined over a number field $k$ and let $G$ be the Galois group of a finite Galois extension $K/k$.
The main aim is to study the action of $G$ on the group $E(K)$.
In this investigation we wish to go deeper than to decompose $E(K)\otimes \CC$ into irreducible representations of $G$, instead we are interested in the integral representation theory of $E(K)$.
Since $\ZZ[G]$-modules can be very complicated even for small groups $G$, we will $p$-adically complete $E(K)$ to get a $\ZZ_p[G]$-module for a prime $p$.
If $p$ does not divide the group order of $G$, we would recover only the structure of $E(K)\otimes \QQ$ and so we will focus on the case when $p$ divides the degree $[K:k]$.

Throughout the paper we will assume that $p$ is an \textbf{odd prime}, that
\textbf{$E(K)$ contains no $p$-torsion elements}
(Assumption~\ref{tors_assu} in Section~\ref{descent_sec}) and that
\textbf{the $p$-primary parts of the Tate-Shafarevich groups $\Sha(E/L)$ are finite}
(Assumption~\ref{sha_assu} in Section~\ref{descent_sec}) for all elliptic curves $E$ and number fields $L$.
The main object of study is the $p$-adic completion $M = E(K)\otimes \ZZ_p$ of $E(K)$ viewed as a $\ZZ_p[G]$-module.

The main aim of our work is to present methods to determine completely $M$ as a $\ZZ_p[G]$-module by using only arithmetic information that is relatively easy to access.
By this we will mean all invariants of $E$ over the global base field $k$ and all local information of $E$ over completions $K_w$ of $K$ at places $w$, but no need for calculations over the global field $K$.
In practice, the knowledge of the type of reduction, the Tamagawa numbers and the number of points in the reduction at a finite number of finite places is all local data that is needed.
In comparison, it is often much harder to calculate explicitly the group $E(K)$: Searching for points over large degree fields $K$ can be very costly and so is bounding the rank by an infinite descent because the class groups involved may be hard to determine.

For instance, if we fix $p=11$ and $G$ to be the dihedral group $D_{11}$, then we can determine the structure of $M$ using only easy information in more than $98.2\%$ of all cases when $E$ has conductor at most $100$ and $K$ is among the ten $D_{11}$-extensions of $\QQ$ with smallest discriminant.
In all these cases, the rank is either $0$ or $11$, but we were never required to find a point in an extension larger than~$\QQ$.

To announce a general theorem in this direction, we need to include a number of conditions, which are mostly generic or could be weakened further.

\begin{thm}\label{ex1_thm}
  Assume that $G$ is isomorphic to the dihedral group $D_p$ with $2p$ elements.
  Let $F$ be the intermediate field with $[F:k]=2$ and let $\twist{E}$ be the quadratic twist of $E$ with respect to~$F$.
  Assume that $\Sha(E/k)$ and $\Sha(\twist{E}/k)$ contain no elements of order $p$, that $K/F$ is everywhere unramified and that all places $v$ such that $p$ divides the Tamagawa number $c_v$ of $E$ or $\twist{E}$ split in $K/F$.
  If $\rk E(k)+\rk \twist{E}(k)\leq 1$, then the $\ZZ_p[G]$-isomorphism class of $M=E(K)\otimes\ZZ_p$ is completely determined by $\rk E(k)$, $\rk \twist{E}(k)$ and local information for $E/K$.
  Furthermore, in all cases $\rk E(K)$ is either $0$, $1$ or $p$.
\end{thm}

One of the reason we restricted here to the dihedral group $D_p$ is that for this group the complete list of indecomposable $\ZZ_p[G]$-modules that are $\ZZ_p$-free only contains six entries.
Even when the conditions of the theorem do not hold, we can very often determine the decomposition of $M$ into a direct sum of these indecomposable modules using only information that is easily accessible.
For $k=\QQ$, we can even show results about the densities with which certain $\ZZ_p[G]$-modules appear (see Theorem~\ref{dihedral_pos_thm} for a precise statement).
See Sections~\ref{ex_cyclic_subsec} and~\ref{ex_dihedral_subsec} for examples illustrating this, and Example~\ref{I_ex} for which we cannot determine $E$ without having to search for points.

For certain groups, including $D_p$, a theorem by Torzewski~\cite{torzewski} lists the invariants we need to know to determine $M$.
(See Proposition~\ref{torzewski2_prop}.)
Our task is to express them in terms of arithmetic invariants of $E$.
The first invariant is the knowledge of $E(K)\otimes \QQ_p$ as a $\QQ_p[G]$-module.
At least conjecturally, we could obtain this from the order of vanishing of $L$-functions of $E$ twisted by the irreducible representations of $G$ via the Birch and Swinnerton-Dyer conjecture.
Instead the other invariants, the so-called Dokchitser regulator constants (see~\eqref{reg_eq} in Section~\ref{metacyc_subsec} for the definition) and the group cohomology $H^1(G,M)$, are not as easy to access.
The main reason that we focus our attention on the dihedral and cyclic case is that for these groups we obtain the best results on how to control the invariants.

As expected the $\ZZ_p[G]$-structure of $M$ is linked to the $p$-primary part  $\Sha(E/K)[p^\infty]$ of the Tate-Shafarevich group.
Here is another general statement that can be deduced from our methods.

\begin{thm}\label{ex2_thm}
  Let $E/\QQ$ be an elliptic curve and let $K/\QQ$ be a cyclic extension of degree $p$.
  Write $u_1$ for the number of primes $\ell$ such that $p$ divides the Tamagawa number $c_\ell$ and $\ell$ is inert in $K/\QQ$.
  Denote by $u_2$ the number of primes $\ell$ that ramify in $K/\QQ$ such that the number of points in the reduction of $E$ over $\ell$ is divisible by $p$.
  If $L(E,\chi,1)\neq 0$ where $\chi$ is a primitive character of $K/\QQ$ and $u_1+u_2 > \rk E(\QQ)$, then the $\FF_p$-dimension of $\Sha(E/K)[p]$ is at least $u_1+u_2-\rk E(\QQ)$.
\end{thm}

We produce this bound by studying the control theorem which links the $p$-primary Selmer group of $E/k$ with the $G$-invariant subgroup of the $p$-primary Selmer group of $E/K$.
The cokernel of the restriction map $\alpha$ between them can be determined completely and we can effectively calculate it using only local information and information of $E$ over $k$.
See Proposition~\ref{cokeralpha_prop} for a precise statement.

One important ingredient for this calculation is to understand the cokernel of the norm map $E(K_w) \to E(k_v)$ where $w$ is a place in $K$ and $v$ the place below~$w$.
This is analogous to the main question in local class field theory and it has been studied before.
In Proposition~\ref{notwild_prop}, we will see that in the case that the ramification index $e_v$ is not divisible by $p$, the cokernel is cyclic of order determined by the Tamagawa number $c_v$ and the residue class degree $f_v$.

This leads to a local question.
The $p$-adic completion $E(K_w)\hatox \ZZ_p$ of $E(K_w)$ is finite unless $K_w$ is a $p$-adic field.
If $K_w$ is a $p$-adic field, then $E(K_w)\otimes \QQ_p$ is isomorphic to $\QQ_p[G_w]$ where $G_w$ is the Galois group of $K_w/k_v$.
In Theorem~\ref{local_thm}, we determine the $\ZZ_p[G_w]$-structure explicitly for all reduction types, when $K_w$ is the unramified cyclic extension of $k_v=\QQ_p$.
Here is a simplified statement not covering all cases.

\begin{thm}\label{ex3_thm}
  Let $E$ be an elliptic curve over $\QQ_p$ with $p\geq 5$ and suppose $K_w/\QQ_p$ is the unramified extension of degree $p$.
  Suppose that the reduction is anything but split multiplicative.
  Then we are in one of the following three cases:
   \begin{itemize}[nosep]
     \item If $E(\QQ_p)$ contains no element of order $p$, then $E(K_w)\hatox \ZZ_p \cong \ZZ_p[G_w]$.
     \item If $E(\QQ_p)$ contains an element of order $p$, but $E(K_w)$ contains no element of order $p^2$, then $E(K_w)\hatox \ZZ_p$ is a non-split extension of $\ZZ_p \oplus \ker\bigl(\ZZ_p[G_w]\to \ZZ_p \bigr)$ by the finite group $\cyclic{p}$ with trivial $G_w$-action.
     \item Otherwise $E(K_w)$ is the direct sum of $\ZZ_p[G_w]$ and a finite group of order $p^2$ with a non-trivial $G_w$-action.
   \end{itemize}
\end{thm}

Underlying to this theorem is the complete classification of all $\ZZ_p[G]$-modules $M$ in the case $G$ is cyclic of order $p$ and $M$ is a finitely generated $\ZZ_p$-module with a cyclic torsion part.
See Theorem~\ref{cyclic_class_thm} for the complete list.

This investigation here grew out of~\cite{bmw1,bmw2} where the $\ZZ_p[G]$-structure of Selmer groups plays an important role in trying to understand the explicit reformulation of the equivariant Birch and Swinnerton-Dyer conjecture under restrictions on the elliptic curve.
One of the motivations of our work is to understand how to remove some of these restrictions; however the present paper does not link to algebraic $L$-values yet.
See~\cite{bley_bsd1,bm_bsd,dokchitser_evans_wiersema}.

Similar methods to the ones used here have been successful in obtaining results on the change of Mordell-Weil groups, Selmer groups and Tate-Shafarevich groups under finite extensions.
For instance, \cite{brau,matsuno,cesnavicius} use versions of the control theorem.
Bartel in~\cite{bartel_large, bartel_psel} used regulator constants to predict growth of the $p$-Selmer group in dihedral extensions.
Our work can be seen as a continuation and generalisation of these methods.

Ouyang and Xie~\cite{ouyang_xie} prove that the size of $\Sha(E/K)$ is unbounded as $K$ varies through the cyclic extensions $K/k$ for a fixed curve $E$ using further methods initiated by Mazur and Rubin as in~\cite{mazur_rubin_finding} and~\cite{mazur_rubin_growth}.
In~\cite{mazur_rubin_ds}, the latter show the following for any elliptic curve $E$ over a number field~$k$: For a positive proportion of primes $p$, for all $n\geq 1$ and all finite sets~$S$ of places in $k$, there are infinitely many cyclic extensions $K/k$ of degree $p^n$ such that all places in $S$ split and $E(K)=E(k)$.
Their emphasis is on constructing extensions $K/k$ given a curve $E/k$, while we fix both and try to determine as much as we can about $E(K)$.

Representation theory, viewing $E(K)\otimes \CC$ as a $\CC[G]$-module, has produced lots of surprising results already.
Much of the work of Tim and Vladimir Dokchitser~\cite{DD09, DD_square} is centred around these questions, especially in connection with parity phenomena.
See~\cite{kellock_dokchitser} for a nice overview with plenty of examples.

In~\cite{bisatt_dokchitser}, the authors use $E(K)\otimes \QQ$ as a $\QQ[G]$-module to make prediction about high order vanishing of certain $L$-functions.
Greenberg has used modular representation theory for Selmer groups in~\cite{greenberg} to obtain results about the growth of the rank in extensions $K/k$, however there $k$ is an infinite extension of $\QQ$.
Finally,~\cite{bley_macias_castillo, macias17} and other work by Macias Castillo and Bley contain the study of $E(K)\otimes \ZZ_p$ as $\ZZ_p[G]$-modules.
Instead the question to determine what $E(K)$ is as a $\ZZ[G]$-module has attracted less attention, likely because it is much harder to say much about it.
However, Theorem~6 in~\cite{dokchitser_evans_wiersema} shows that the arithmetic of $L$-values should predict interesting results in this direction.

The structure of this article is as follows.
In Section~\ref{local_norm_sec} we will investigate the cokernel of the norm map for a local extension.
Part of these results are well-known, but we try to be as general as possible.
Then Section~\ref{rep_sec} is devoted to gathering results on the integral representation theory for a cyclic group of order $p$, where we allow non-trivial torsion, and certain groups that are extensions of cyclic groups of order coprime to~$p$ by cyclic groups of order~$p$.
This is then used in Section~\ref{local_mw_sec} where we determine the local group of points as a $\ZZ_p[G_w]$-module in the case of the unramified extension of degree $p$.
General results on the control theorem in a general global extension are presented in Section~\ref{descent_sec}.
Finally, they are applied to global extensions with cyclic or dihedral Galois groups in Section~\ref{galmod_sec}.
This section also includes a list of examples, which illustrate how to use the general method to determine the $\ZZ_p[G]$-module structure.

\subsection*{Acknowledgments}
We wish to thank the anonymous referee for the careful reading of the article which helped us to improve the article.
We are also indebted to Lewis Matthews for valuable feedback.

\section{Notations}

This section is not meant to be read as such, but it can help to find the definitions of the notations commonly used.

Throughout, \hl{$p$} is an odd prime and \hl{$\ZZ_p$} denotes the ring of $p$-adic integers.

In general, for an abelian group $Z$, we denote the projective limit of $Z/p^nZ$ by \hl{$Z\hatox \ZZ_p$}, which we call the \hl{$p$-adic completion} of $Z$.
If $Z$ is finitely generated this coincides with $Z\otimes_{\ZZ} \ZZ_p$ and, if $Z$ is finite, it is isomorphic to the \hl{$p$-primary} torsion subgroup \hl{$Z[p^\infty]$}.
If $Z$ is a discrete abelian $p$-primary group, or a compact $\ZZ_p$-module, then \hl{$Z^{\vee}$} denotes the \hl{Pontryagin dual} $\Hom(Z,\QZp)$.
At times, we will omit the ring $\ZZ$ or $\ZZ_p$ over which we take tensor products: specifically, if we tensor a Mordell-Weil group of an elliptic curve over a number field, we will always mean $\otimes_{\ZZ}$, instead if $M$ is a $\ZZ_p$-module then $M\otimes \QQ_p$ and $M\otimes \QZp$ mean $\otimes_{\ZZ_p}$.

Throughout the paper \hl{$G$} will stand for a finite group.
For a $\ZZ_p[G]$-module $M$, we will write \hl{$M_t$} for the torsion subgroup of $M$ and \hl{$M_f$} for the quotient of $M$ by $M_t$.
The action is always from the left, even if we tend to write $M\otimes \QQ_p$ for the $\QQ_p[G]$-module obtained by extending the scalars.
In Section~\ref{saturation_subsec}, we will define the \hl{saturation index $\iota(M)$}.
The \hl{Dokchitser regulator constants $\mathcal{C}_{\Theta}(M)$} and their valuation \hl{$s_{\Theta}(M)$} are defined in Section~\ref{metacyc_subsec}.
The groups \hl{$H^i(G,M)$} refer to the usual group cohomology with $\hl{M^G}=H^0(G,M)$ and \hl{$\hat H^i(G,M)$} is the modified version by Tate.

The symbols \hl{$\ZZ_p$}, \hl{$\twist{\ZZ}_p$}, \hl{$\ZZ_p\{i\}$}, \hl{$A$}, \hl{$\twist{A}$}, \hl{$A\{i\}$}, \hl{$B$}, \hl{$\twist{B}$} and \hl{$B\{i\}$} denote indecomposable $\ZZ_p[G]$-modules for the case when $G$ is a cyclic group of order $p$ as in Section~\ref{cyclicrep_subsec}, a dihedral group $D_p$ or one of the metacyclic groups as in Section~\ref{metacyc_subsec}.
Furthermore by \hl{$\bigl\{\cyclic{p^i}\bigm\vert \ZZ_p\bigr\}$}, \hl{$\bigl\{\cyclic{p^i}\bigm\vert A \bigr\}$} and, \hl{$\bigl\{\cyclic{p^i}\bigm\vert \ZZ_p\oplus A\bigr\}$} we denote certain non-split extensions in the case of cyclic groups as explained in Theorem~\ref{cyclic_class_thm}.
Proposition~\ref{cyclic_fp_prop} and Lemma~\ref{cyclic_fin_lem} contain the definitions of the finite modules \hl{$J_i$} and \hl{$F_{i,\bar{w}}$}.

The symbols \hl{$k$} and \hl{$K$} will stand for fields such that $K/k$ is an Galois extension of group~$G$.
In Sections~\ref{local_norm_sec} and~\ref{local_mw_sec} they are local fields, while they are number fields in the Sections~\ref{descent_sec} and~\ref{galmod_sec}.
For a place $v$ in $k$, the completion is \hl{$k_v$} and \hl{$\FF_v$} is the residue field.
The discriminant is \hl{$\Delta_k$}.
The letters \hl{$e=e_v$} and \hl{$f=f_v$} are the ramification index and residue class degree at a place $v$.

The letter \hl{$E$} will stand for an elliptic curve defined over $k$, while \hl{$\mathcal{E}$}, \hl{$\mathcal{E}^0_v$}, and \hl{$\Phi_v$} relate to the Néron model as explained at the start of Section~\ref{local_norm_sec}.
If the reduction is good, we will use \hl{$\tilde{E}(\FF_v)$} for the group of points in the reduction.
The  Tamagawa number of $E$ at the finite place $v$ is denoted by  \hl{$c_v$}.
The modified product  \hl{$C(E/k)$} of the Tamagawa number is defined in Section~\ref{galmod_sec}.

For a finite place $v$, the group $\hl{D=D_v}= \hat H^0\bigl(G, E(K_w)\hatox \ZZ_p\bigr)$ is investigated in Section~\ref{local_norm_sec}.
The \hl{$p$-primary Selmer group} \hl{$\Sel_k$} of $E$ appears first in Section~\ref{descent_sec}.
To shorten the notation, we will write \hl{$\Sha_k=\Sha(E/k)[p^\infty]$} for the $p$-primary part of the \hl{Tate-Shafarevich group}.
In Section~\ref{descent_sec}, we will encounter the maps \hl{$\alpha$}, \hl{$\beta$}, \hl{$\gamma$}, \hl{$\delta$}, \hl{$\varepsilon$} and \hl{$\eta$}.
The \hl{capitulation} is the kernel \hl{$C_{K/k} = \ker (\eta)$} and \hl{$D_{K/k}$} stands for the sum of $D_v$ over all places in the set \hl{$S$} containing all places of bad reduction for $E$, all places ramified in $K/k$ and all infinite places.

In Section~\ref{galmod_sec}, by \hl{$\twist{E}$} we denote a quadratic twist of $E$.
While \hl{$r_F$} is the rank of $E(F)$ for a field $F$, the rank of $\twist{E}(F)$ is \hl{$\twist{r}_F$}.

\section{The local norm}\label{local_norm_sec}

The aim of this section is to study the cokernel of the norm map on an elliptic curve under a finite extension of local fields.
The analogous question, which is central in local class field theory, concerns the cokernel of the norm map on units.
However the situation is more complicated here and it will lead us to treat different cases separately.
While this has been studied partially in many situations, we try to be as general as possible.

Let $k$ be a local field with valuation $v$ and let $K$ be a finite Galois extension with valuation $w$.
By $\FF_v$ and $\FF_w$ we will denote their residue fields, by $\mathcal{O}_v$ and $\mathcal{O}_w$ their rings of integers and by $\mathfrak{m}_v$ and $\mathfrak{m}_w$ their maximal ideals.
Let $G$ be the Galois group of $K/k$.
The ramification index is $e$ and $f$ stands for the residue class degree of the extension $K/k$.
The degree $[K:k]=e\cdot f$ is denoted by $n$.
When we write $v\mid p$, we mean that $k$ is a finite extension of $\QQ_p$.
The group we wish to determine is
\[
  D = \hat{H}^0\bigl(G, E(K)\hatox\ZZ_p\bigr) = \coker\bigl( \No\colon E(K)\hatox \ZZ_p\to E(k)\hatox \ZZ_p\bigr)
\]
where $\hat H$ denotes Tate's modification of group cohomology and $E(K)\hatox \ZZ_p$  is the $p$-adic completion $\varprojlim E(K)/p^i E(K)$ of $E(K)$.

Let $\mathcal{E}_v$ be the Néron model of $E$ over the ring of integers $\mathcal{O}_v$ of $k$.
Write $\mathcal{E}^0_v/\mathcal{O}_v$ for the connected component of the identity, $\Ens_v^0/\FF_v$ for its special fibre, and let $\Phi_v/\FF_v$ be the group of components of the special fibre.
If we have good reduction, we will simply write $\tilde E$ for $\mathcal{E}^0_v$ as there is no danger of confusion.
The Tamagawa number of $E$ over $k$ is denoted by $c_v = \lvert \Phi_v(\FF_v)\rvert$.
We use similar notation for $E/K$ with $v$ replaced by $w$.

\begin{lem}\label{reduction_lem}
  If $v\nmid p$, then the Néron model $\mathcal{E}_w$ over $\mathcal{O}_w$ can be chosen such that we have an exact sequence of finite $\ZZ_p[G]$-module
  \[
    \xymatrix@1{
     0\ar[r] & \Ens^0_w(\FF_w)[p^\infty] \ar[r] & E(K)\hatox\ZZ_p\ar[r] & \Phi_w(\FF_w)[p^\infty]\ar[r] & 0.
    }
  \]
\end{lem}
Note this does not necessarily mean that the $G$-fixed part of the outer terms of this sequence are the corresponding groups for $k$ as the Néron model may change in the extension.
\begin{proof}
  We claim that, even if the Néron model changes, we may choose the model over $\mathcal{O}_w$ such that the subgroup $\mathcal{E}_w^0(\mathcal{O}_w)$ is a $G$-submodule.
  To prove this claim, we may assume that $E$ has bad reduction over $K$ and hence over $k$, too.
  We may translate a chosen equation over $\mathcal{O}_v$ to obtain a first Weierstrass model whose singular point over $\FF_v$ is $(0,0)$.
  The Néron model over $\mathcal{O}_w$ constructed starting from this equation satisfies the claim.

  Since the kernel of $\mathcal{E}^0_w(\mathcal{O}_w) \to \Ens^0_w(\FF_w)$ is divisible by $p$, we have $\mathcal{E}^0_w(\mathcal{O}_w)\hatox \ZZ_p \cong \Ens^0_w(\FF_w)\hatox\ZZ_p = \Ens^0_w(\FF_w)[p^\infty]$.
  The projective limit over $m$ of the exact sequence
  \[
    \xymatrix@1{ \Phi_w(\FF_w)[p^m]\ar[r] & \mathcal{E}^0_w(\mathcal{O}_w)/p^m \ar[r] & E(K)/p^m \ar[r] & \Phi_w(\FF_w)/p^m\ar[r] & 0 }
  \]
  stays exact and the first term will vanish as $\Phi_w$ is finite.
\end{proof}

\begin{lem}\label{mult_lem}
  Suppose $E$ has split multiplicative reduction over $k$.
  Let $q \in\mathcal{O}_v$ be the Tate parameter such that $E(k) \cong  k^\times/q^{\ZZ}$.
  Then $D$ is the quotient of the $p$-primary part of $k^{\times}/\No(K^\times)$ by the subgroup generated by the image of~$q$.
\end{lem}
Write $G^{p\text{-ab}}$ for the Galois group of the maximal abelian $p$-extension within $K/k$.
If $\operatorname{rec}\colon k^{\times}\to G^{p\text{-ab}}$ is the reciprocity map, then $D$ is isomorphic to $G^{p\text{-ab}}/\langle\operatorname{rec}(q)\rangle$.
\begin{proof}
  First recall that the reduction type of $E$ is still split multiplicative over $K$ with the same parameter $q$.
  Since the torsion subgroup of $E(K)$ is finite, we get a diagram with exact rows
  \[
    \xymatrix@R-1ex{
      0\ar[r] & q^{\ZZ_p}\ar[r]\ar[d] & K^{\times}\hatox\ZZ_p\ar[r]\ar[d]^{\No} & E(K)\hatox\ZZ_p\ar[r] \ar[d]^{\No} & 0 \\
      0\ar[r] & q^{\ZZ_p}\ar[r] & k^{\times}\hatox\ZZ_p\ar[r] & E(k)\hatox\ZZ_p\ar[r] & 0. \\
    }
  \]
  This shows that $D$, the cokernel on the right, is the quotient of the cokernel in the middle by the subgroup generated by the image of~$q$.
\end{proof}

\begin{prop}
  Suppose $E$ has split multiplicative reduction over $k$ and that $v\nmid p$.
  \begin{itemize}
    \item If $p\mid c_v$ and $p\mid f$, then $D$ is non-trivial.
    \item If $p\nmid \gcd(f,c_v)$ and $p\nmid \gcd\bigl(e,\lvert\FF_v^\times\rvert\bigr)$, then $D$ is trivial.
  \end{itemize}
\end{prop}

Note this leaves the case when $p$ divides $e$ and $\lvert \FF_v^\times\rvert$, but does not divide $c_v$ or~$f$.
In that last case, it could be trivial or non-trivial, which can only be determined by a finer analysis.

\begin{proof}
  Note that the reduction is still split multiplicative over $K$.
  We rewrite the sequence from Lemma~\ref{reduction_lem} in this situation and compare it to the same sequence but over $k$.
  This gives us the following commutative diagram:
  \[
    \xymatrix{%
    0\ar[r] & \FF_w^{\times}[p^\infty] \ar[r]\ar[d]^{[e]\cdot \No_{\FF_w/\FF_v}} &
    E(K)\hatox \ZZ_p \ar[r] \ar[d]^{\No} & \cyclic{ec_v}\ar[r]\ar[d] & 0 \\
    0\ar[r] & \FF_v^{\times}[p^\infty] \ar[r] &
    E(k)\hatox \ZZ_p \ar[r] & \cyclic{c_v}\ar[r] & 0. \\
    }
  \]
  The vertical map on the right sends $1+ec_v\ZZ$ to $f+c_v\ZZ$.
  We deduce the exact sequence
  \[ \xymatrix@1{
    \cyclic{e\,(c_v, f)} \ar[r] & \FF_v^{\times}[p^{\infty}]/e \ar[r] & D \ar[r] & \cyclic{(c_v,f)}\ar[r] & 0
    }
  \]
  where the second term is the quotient of the $p$-primary component of $\FF_v^{\times}$ by its $e$-th powers.
  If $p\mid \gcd(c_v,f)$, then $D$ is non-trivial.
  If both terms next to $D$ are trivial, i.e., when $p\nmid\gcd(c_v,f)$ and $p\nmid \gcd(e, \lvert \FF_v^\times\rvert)$, then $D$ is trivial.
\end{proof}

\subsection{Unramified places}

\begin{lem}\label{unram_lem}
  Suppose $K/k$ is unramified.
  Then $D \cong \hat H^0\bigl( G, \Phi_v(\FF_w)[p^\infty]\bigr)$.
\end{lem}
\begin{proof}
  Since we assume that the extension is unramified, the Néron model does not change: $\mathcal{E}_w = \mathcal{E}_v\times \mathcal{O}_w$ and in particular $\Phi_w = \Phi_v\times \FF_w$.
  By the three propositions in Section~4.1 of~\cite{ellerbrock_nickel}, the Galois module $\mathcal{E}^0_w(\mathcal{O}_w)$ is cohomologically trivial for any unramified extension.
  From the fact that $\Phi_w(\FF_w)=\Phi_v(\FF_w)$ is finite, we obtain the exact sequence
  \[
    \xymatrix@1{ 0\ar[r] & \mathcal{E}_w^0(\mathcal{O}_w)\hatox \ZZ_p \ar[r] & E(K)\hatox\ZZ_p\ar[r] & \Phi_v(\FF_w)[p^\infty]\ar[r] & 0.
    }
  \]
  Since the first term is cohomologically trivial, we get $D \cong \hat H^0\bigl(G, E(K)\hatox\ZZ_p\bigr)$ is isomorphic to
  $\hat H^0\bigl( G, \Phi_v(\FF_w)[p^\infty]\bigr)$.
\end{proof}

For any integer $m$, we will denote the highest power of $p$ dividing $m$ by $\gcd(m,p^\infty)$.

\begin{lem}\label{unram_dv_lem}
  Assume that $K/k$ is unramified.
  Then $D$ is a cyclic group of order $c' = \gcd\bigl(c_v,n,p^\infty\bigr)$.

  In particular, $D$ is trivial except possibly in the following two cases:
  \begin{itemize}[nosep]
  \item $E$ has split multiplicative reduction over $k$ of Kodaira type I${}_n$ with $p\mid n$ or
  \item $p=3$ and the special fibre of $\mathcal{E}_v$ is of Kodaira type IV or IV${}^*$.
  \end{itemize}
\end{lem}

\begin{proof}
  Since $K/k$ is unramified $\Phi_w = \Phi_v \times \FF_w$ and $D$ is a quotient of $\bigl(\Phi_v(\FF_w)[p^\infty]\bigr)^{G} = \Phi_v(\FF_v)[p^\infty]$ by Lemma~\ref{unram_lem}.
  Therefore, if $p$ does not divide $c_v$ then $D$ is trivial.
  From the assumptions that $p\neq 2$, the classification of bad fibres of elliptic curves implies that $\Phi_v(\FF_w)$ may contain a $p$-torsion element only if the reduction is split multiplicative or if $p=3$ and the type is IV or IV${}^*$.

  Assume that $E$ has split multiplicative reduction in which case we may use Lemma~\ref{mult_lem}.
  Since the extension is unramified, the valuation $v$ induces an isomorphism $k^\times /\No(K^\times)$ to $\cyclic{n}$.
  As the valuation of the Tate parameter $q$ is equal to the Tamagawa number $c_v$, we find that $D$ is indeed cyclic of order $c'$.

  If $p=3$ and the fibre is of type IV or IV${}^*$ with $p=c_v$ then $\Phi_v$ is the constant group scheme~$\cyclic{3}$.
  Therefore $D\cong \hat H^0\bigl(G, \cyclic{3}\bigr)$ is the cokernel of multiplication by $n$ on $\cyclic{3}$.
  So once again $D$ is cyclic of order $\gcd\bigl(c_v,n,p^\infty\bigr)$.
\end{proof}

 If $p$ were allowed to be $2$, then we could in the same fashion go through all Kodaira types in Tate's algorithm and determine explicitly the group $D$ from the type and the degree $n$.

\subsection{Totally ramified places}\label{tot_ram_subsec}

In this subsection we assume $K/k$ is totally ramified.
We begin by considering the case where $v$ does not divide $p$.

\begin{prop}\label{non_p_totram_prop}
  Suppose that $v\nmid p$.
  \begin{itemize}
    \item If the reduction of $E$ is good then, then $D \cong Z /n Z$ where $Z$ is the group $\tilde E(\FF_v)[p^\infty]$.
    \item Suppose $E$ has split multiplicative reduction with Tamagawa number $c_v$ and Tate parameter~$q$.
    Write $q = u\cdot \No(\pi_w)^{c_v}$ for a choice of a uniformiser $\pi_w$ of $K$ and a unit $u\in \mathcal{O}_v^\times$.
    Then $D$ is isomorphic to the quotient of the $p$-primary part of $\FF_v^\times$ by its subgroup generated by the image of~$u$ and by all $n$-th powers.
    \item If $E$ has non-split multiplicative reduction over $k$, then $D$ is cyclic of order $\gcd\bigl(n,\lvert \FF_v \rvert +1, p^\infty\bigr)$.
    \item If $E$ has additive reduction, then $D$ is cyclic of order $\gcd\bigl(c_v,n,p^{\infty}\bigr)$.
    \end{itemize}
\end{prop}

\begin{proof}
  Suppose first that the reduction is good over $k$.
  Then the reduction is also good over $\mathcal{O}_w$ and Lemma~\ref{reduction_lem} tells us that $D\cong\hat H^0\bigl(G, \tilde E(\FF_w)[p^\infty]\bigr)$.
  As the extension is totally ramified the action of $G$ is trivial on $Z=\tilde E(\FF_w)[p^\infty]$.
  We conclude that $D$ is isomorphic to $Z/n\,Z$.

  If the reduction is split multiplicative over $k$ we use Lemma~\ref{mult_lem}.
  Since the extension is totally ramified and $v\nmid p$, the $p$-primary part of the group $k^\times/\No(K^\times)$ identifies with the $p$-primary part of the quotient of $\FF_v^\times$ by its $n$-th powers.
  The group generated by the image of~$q$ under this identification is the one generated by the image of~$u$.

  Next, we treat the case when the reduction is additive with $p\nmid c_v$.
  Then $E(k)\hatox\ZZ_p \cong \Ens_v^0(\FF_v)[p^\infty]$ by Lemma~\ref{reduction_lem}; this group is trivial as $v\nmid p$ and therefore $D$ is also trivial.

  The same argument also works for non-split multiplicative reduction as the reduction will still be non-split over $K$, except that $\Ens_v^0(\FF_v)$ is now a cyclic group of order $\lvert \FF_v\rvert+1$ with the trivial action by $G$ on it.

  Finally, we are left with the case of reduction of type IV or IV${}^*$ and $c_v=p=3$.
  If the reduction is still of the same type over $K$, then $E(K)\hatox \ZZ_p$ is isomorphic to $\Phi_w(\FF_w)\cong \cyclic{3}$ with trivial action by $G$.
  Hence in this case $D$ is trivial unless $3\mid n$ in which case it is cyclic of order $3=c_v$.

  Instead if the reduction type changes, the Tamagawa number $c_w$ must be coprime to~$3$.
  Then $E(K)\hatox\ZZ_p$ is isomorphic to $\Ens^0_w(\FF_w)[p^\infty]$ and has trivial action by $G$.
  But the $G$-fixed part is $E(k)\hatox\ZZ_p$ which is cyclic of order $3$.
  Yet again we can draw the same conclusion as before.
\end{proof}

Now we consider one case when $v$ divides $p$.
First we give a general lemma.

\begin{lem}\label{coinf_cores_lem}
  Let $X$ be a $G$-module and let $H\leq G$ be a normal subgroup.
  Then the following is an exact sequence:
  \[
    \xymatrix@1{\hat H^{-1}(G,X)\ar[r] & \hat H^{-1}(G/H, X^H)\ar[r] & \hat H^0(H,X)_{G/H}\ar[r] & \hat H^0(G,X)\ar[r] & \hat H^0(G/H,X^H)\ar[r] & 0 }
  \]
  In particular, if $X^H$ is cohomologically trivial as a $G/H$-module, then $\hat H^0(G,X)$ is isomorphic to $\hat H^0(H,X)_{G/H}$.
\end{lem}
\begin{proof}
  The norm map with respect to $G$ is equal to the composition
  \[ X_G=(X_H)_{G/H} \to (X^H)_{G/H} \to (X^H)^{G/H} = X^G\]
  where the first map $\rho$ is induced by the norm map for $H$ and the second is the norm map for $G/H$.
  The result can now be deduced from the kernel-cokernel sequence, which gives the desired exact sequence except that the term $\hat H^0(H,X)_{G/H}$ is replaced by $\coker(\rho)$.
  But these two groups are equal as taking $G/H$-coinvariants is right exact.
\end{proof}

\begin{prop}\label{p_ram_prop}
  Suppose $v\mid p$.
  Assume that $K/k$ is totally ramified and that $E$ has good ordinary reduction over $k$.
  Let $Z=\tilde E(\FF_v)[p^\infty]$ and let $Y=G^{p-\text{ab}}$ be the maximal abelian $p$-primary quotient of $G$.
  Then there is an exact sequence
  \[
    \xymatrix@1{
      Z\bigl[ n \bigr] \ar[r] & Y /(1 - \alpha)Y \ar[r] & D \ar[r] & Z / n\, Z \ar[r] & 0,
    }
  \]
  where $\alpha \in \ZZ_p^\times$ is the unit root of the characteristic polynomial $X^2 - a \,X + \#\FF_v = 0$  of Frobenius with $a = \#\FF_v + 1 - \# \tilde{E}(\FF_v)$.
\end{prop}

\begin{proof}
  Write $\hat{E}$ for the formal group associated to a minimal Weierstrass equation of $E$.
  Consider the exact sequence

  \[\mathclap{\xymatrix@1{
      \hat H^{-1}\bigl(G, \tilde E(\FF_w)[p^\infty]\bigr)\ar[r] & \hat H^0\bigl(G,\hat E(\mathfrak{m}_w)\bigr)\ar[r] & \hat H^0\bigl(G,E(K)\hatox\ZZ_p\bigr)\ar[r] & \hat H^0\bigl(G, \tilde E(\FF_w)[p^\infty]\bigr)
    }}
  \]
  where the very last map is surjective since the map $E(k)\hatox\ZZ_p\to \tilde E(\FF_v)[p^\infty]$ is surjective.
  As the extension is totally ramified, the two extremal non-trivial terms in the sequence are isomorphic to the kernel and cokernel of multiplication by $n$ on the group $Z$.

  Let $H$ be the $p$-Sylow subgroup of $G$, i.e, the wild inertia subgroup, and let $L = K^H$ with valuation $w'$.
  Since $G/H$ is of order coprime to $p$, the cohomology $\hat{H}^i\bigl(G/H, \hat{E}(\mathfrak{m}_{w'})\bigr)$ must vanish (see VI.5.5 in~\cite{brown}).
  We apply the previous lemma with $X= \hat{E}(\mathfrak{m}_w)$ to obtain that $\hat{H}^0\bigl(G,\hat{E}(\mathfrak{m}_w)\bigr)$ is isomorphic to the largest quotient of $\hat{H}^0\bigl(H, \hat{E}(\mathfrak{m}_w)\bigr)$ on which $G/H$ acts trivially.

  Theorem~1 in~\cite{lubinrosen} shows that $\hat{H}^0\bigl( H, \hat{E}(\mathfrak{m}_w)\bigr)$ is isomorphic to $H^{\text{ab}}/(1-\alpha)H^{\text{ab}}$ where $H^{\text{ab}}$ is the abelianisation of $H$ viewed as a $\ZZ_p$-module.
  See also Appendix~B in~\cite{mazur_rubin_finding}, which includes a reference that $\alpha$ is the twisting matrix in~\cite{lubinrosen}.
  This isomorphism is $G/H$-equivariant by inspecting its construction in~\cite{lubinrosen} using the formula in Proposition~IV.9 in~\cite{serre_cl} with $i=1$.
  Therefore $\hat{H}^0\bigl(G,\hat{E}(\mathfrak{m}_w)\bigr)$ is isomorphic to the largest quotient of  $H^{\text{ab}}/(1-\alpha)H^{\text{ab}}$ on which $G/H$ acts trivially.
  The largest quotient of $H^{\text{ab}}$ on which $G/H$ acts trivially by conjugation is the Galois group of the largest extension $L'$ of $L$ which is abelian, $p$-primary and such that the Galois group of $L'/k$ is the direct product of $H$ and $\Gal(L'/L)$.
  It is therefore also the Galois group of the largest extension of~$k$, which is abelian and $p$-primary, i.e., it is $Y$.
  As taking $G/H$-coinvariants is right exact, the $G/H$-coinvariants of $H^{\text{ab}}/(1-\alpha)H^{\text{ab}}$ is isomorphic to $Y/(1-\alpha)Y$.
\end{proof}

Recall that we say that $E$ has anomalous reduction at a place of good reduction over $k$ if $p$ divides $\# \tilde E(\FF_v)$.

\begin{cor}\label{p_ram_cor}
  Suppose $v\mid p$.
  Assume that $K/k$ is totally and wildly ramified and that $E$ has good ordinary reduction over $k$.
  Then $D$ is non-trivial if and only if $E$ has anomalous reduction.
\end{cor}

\begin{proof}
  The group $Z$ in Proposition~\ref{p_ram_prop} is non-trivial if and only if $E$ has anomalous reduction.
  In that case $D$ is non-trivial as $n$ is divisible by $p$.

  Instead assume that the reduction is not anomalous.
  Then $1-\alpha$ divides the evaluation of the characteristic polynomial at $X=1$, which equals $\# \tilde{E}(\FF_v)$.
  Therefore $1-\alpha\in\ZZ_p^{\times}$ and therefore $D$ is trivial.
\end{proof}

\subsection{The general case}

We drop any assumption on the local extension $K/k$.
Recall that $e$ is the ramification index, which is the order of the inertia subgroup $I$.

\begin{prop}\label{notwild_prop}
  If $p\nmid e$, then $D$ is cyclic of order $c' = \gcd(c_v,n,p^\infty)$.
\end{prop}
\begin{proof}
  Apply Lemma~\ref{coinf_cores_lem} with $X=E(K)\hatox \ZZ_p$ and $H=I$.
  The assumption that $p\nmid e$ implies that $\hat H^0(I,X)=0$.
  Therefore $D$ is reduced to the computation of the same group for the unramified extension $K^{I}/k$ and that was done in Lemma~\ref{unram_dv_lem}.
\end{proof}

More generally, this proof shows that $D$ has a cyclic quotient of order $c'=\gcd(c_v,f,p^{\infty})$ where $f$ is the residue class degree, for all local extensions.

\begin{prop}\label{d_coinv_prop}
  Suppose that $p\nmid c_v$.
  If $p=3$, assume further that $f$ is odd or that the reduction type is not IV or IV${}^*$.
  If the reduction is non-split multiplicative of type I${}_n$, assume that $f$ is odd or $p\nmid n$.
  Then $D\cong\hat H^0\bigl(I, E(K)\hatox\ZZ_p\bigr)_{G/I}$.
\end{prop}

\begin{proof}
   Let $L$ be the subextension fixed by $I$ with valuation $w'$.
   The assumptions are made to assure that the Tamagawa number of $E$ over $L$ is still not divisible by $p$:
   If the Tamagawa number $c_{w'}$ over $L$ is divisible by $p$, but not $c_v$, then we must either be in the case IV or IV${}^*$ and $p=3$ or in the case I${}_n$ with $p\mid n$.
   Further the group $\Phi_v$ acquires new points of order $p$ in either cases, only if $L/k$ is of even degree.

   We are going to use Lemma~\ref{coinf_cores_lem} with $X = E(K)\hatox \ZZ_p$ and $H=I$.
   Since $L/k$ is unramified, the Néron model $\mathcal{E}_v$ does not change under the extension.
   Therefore $X^{I} = E(L)\hatox \ZZ_p = \mathcal{E}^0_v(\mathcal{O}_{w'})\hatox \ZZ_p$ is a $G/I$-module that is cohomologically trivial, again by~\cite{ellerbrock_nickel}.
   This concludes the proof.
\end{proof}

In practice this means that in these cases the calculation of $D$ reduces to the totally ramified case treated in Section~\ref{tot_ram_subsec}.
Here is one important example.

\begin{prop}\label{good_not_p_prop}
  Suppose $E$ has good reduction and $v\nmid p$.
  Then $D \cong Z/eZ$ where $Z=\tilde E(\FF_v)[p^\infty]$. 
\end{prop}
\begin{proof}
  From the proof of Proposition~\ref{non_p_totram_prop}, we see that $\hat H^0\bigl(I, E(K)\hatox\ZZ_p\bigr)$ is isomorphic to the group $\tilde E(\FF_w)[p^\infty]/e\tilde E(\FF_w)[p^\infty]$ as a $G/I$-module.
  By Proposition~\ref{d_coinv_prop}, we want to calculate its $G/I$-coinvariants.
  Lang's theorem~\cite{lang} implies that $H^1\bigl(\FF_w/\FF_v,\tilde E(\FF_w)\bigr)=0$ and that the norm $\tilde E(\FF_w)\to \tilde E(\FF_v)$ is surjective as the Herbrand quotient of the finite $\tilde E(\FF_w)$ is~$1$.
  Therefore the $G/I$-coinvariant space of $\tilde E(\FF_w)[p^\infty]$ is $Z=\tilde E(\FF_v)[p^\infty]$.
  We conclude that the $G/I$-coinvariants of $\tilde E(\FF_w)[p^\infty]/e\tilde E(\FF_w)[p^\infty]$ is $Z/eZ$ because taking $G/I$-coinvariants is right exact.
\end{proof}

While the results in this section do not cover all cases, they do cover a lot and in the remaining ones one can often use the same methods to reduce it to a simple calculation.
There is one major exception to this, that is when we have wild ramification, but the reduction is not good ordinary.
Of course, in practice it is always possible to calculate explicitly a generating set for $E(K)$ and to determine $D$ from its definition; this will involve Tate's algorithm and Hensel's lifting over $\mathcal{O}_w$.

\section{Representation theory}\label{rep_sec}

In his section, we will gather the representation theoretic results used later in the case that the Galois group is either cyclic or dihedral.
We should emphasis that this is integral representation theory, in that we deal with $\ZZ_p[G]$-modules.

By the Krull-Schmidt-Azumaya Theorem (see Theorem~6.12 in~\cite{cr1}) any finitely generated $\ZZ_p[G]$-module decomposes in a unique way into a direct sum of indecomposable $\ZZ_p[G]$-modules.
For a finite group $G$, we will say that $M$ is a $\ZZ_p[G]$-lattice, or simply $G$-lattice, if it is a finitely generated $\ZZ_p[G]$-module that is free as a $\ZZ_p$-module.
As all modules will be finitely generated, we say $M$ is a finite module if it has finitely many elements, and no confusion should arise.

\subsection{Cyclic group of order \texorpdfstring{$p$}{\textit{p}}}\label{cyclicrep_subsec}

Let $G$ be the cyclic group of order $p$.
We will write $\tau$ for a generator of $G$.

This is of course the easiest group of interest and, not surprisingly, we can give clear classification results for $\ZZ_p[G]$-modules.
First, there are only two irreducible $\QQ_p[G]$-modules, namely $\QQ_p$ and the kernel of $\No\colon\QQ_p[G]\to \QQ_p$.
The classification of indecomposable $\ZZ_p[G]$-lattices is also well-known (see Section~34B of~\cite{cr1}, though this goes back to~\cite{diederichsen} and~\cite{reinercyclic}).

\begin{prop}\label{cyclic_lattices_prop}
  There are exactly three isomorphism classes of indecomposable $\ZZ_p[G]$-lattices, namely the trivial lattice $\ZZ_p$, the free module $\ZZ_p[G]$ and the augmentation kernel $A$ of $\No\colon \ZZ_p[G]\to \ZZ_p$.
\end{prop}
We may view $A=(\tau-1)\ZZ_p[G]$ also as the ring $\ZZ_p[\zeta]$ with $\zeta$ a primitive $p$-th root of unity and the action of $\tau\in G$ given by multiplication with $\zeta$.
It is important to note that $\ZZ_p[G]$ and $\ZZ_p\oplus A$ are the only $G$-lattices $M$ with $M\otimes \QQ_p\cong \QQ_p[G]$.
In this simple case the concepts of free and projective coincide.

\begin{prop}\label{cyclic_fp_prop}
  For $1\leq i\leq p$, let $J_i$ be the $i$-dimensional $\FF_p$-vector spaces and the action by $\tau$ written as a unique Jordan block with $1$ on the diagonal.
  Then these $J_i$ are the only indecomposable finitely generated $\FF_p[G]$-modules.
\end{prop}

\begin{proof}
  The action of $G$ on a finite dimensional $\FF_p$-vector space is described by the matrix of the action by the generator $\tau$.
  This matrix can be put into Jordan normal form.
  The eigenvalues of a matrix of order $p$ must be $1$ and the module is indecomposable if there is only one block.
  Such a single Jordan block has order $p$ if and only if its dimension is at most~$p$.
\end{proof}

For any $i\geq 1$, by $\cyclic{p^i}$ we will mean the cyclic group of order $p^i$ with trivial action by $G$.
For $i\geq 2$ and $\bar w \in \FF_p^{\times}$, set $w = 1+p^{i-1}\bar{w}\in \cyclic{p^i}$.
The group homomorphism $G\to \bigl(\cyclic{p^i}\bigr)^{\times}$ sending $\tau$ to $w$ defines a non-trivial action of $G$ on the cyclic group of order $p^i$.
We will write $F_{i,\bar{w}}$ for this $G$-module.
These $p-1$ non-trivial $\cyclic{p^i}[G]$-modules $F_{i,\bar{w}}$ are pairwise non-isomorphic.

\begin{lem}\label{cyclic_fin_lem}
  Let $i\geq j\geq 0$.
  If $M$ is a finite $\ZZ_p[G]$-module, which is cyclic of order $p^i$ such that $M^G$ is cyclic of order $p^j$, then $M$ is either $\cyclic{p^i}$ and $i=j$ or $M$ is $F_{i,\bar{w}}$ with $i = j+1\geq 2$ and $\bar{w}\in\FF_p^{\times}$.
\end{lem}
\begin{proof}
  If $j<i$, the action is non-trivial and we must have $M\cong F_{i,\bar{w}}$.
  Then $M^G = pM$ shows that $i=j+1$.
\end{proof}

We are now interested in $\ZZ_p[G]$-modules $M$ whose torsion part $M_t$ is cyclic.
We will use the notation $\{X\vert Y\}$ representing a non-split extension $\xymatrix@1{0\ar[r]&X\ar[r] &M\ar[r] & Y \ar[r] & 0}$.
The ones appearing in the following theorem will be constructed explicitly in its proof.

\begin{thm}\label{cyclic_class_thm}
  Let $M$ be a finitely generated $\ZZ_p[G]$-module and suppose that its $\ZZ_p$-torsion subgroup $M_t$ is cyclic.
  Then $M$ is a direct sum of some of the following indecomposable modules:
  \begin{center}\renewcommand{\arraystretch}{1.3}
    \begin{tabular}{cccl}
    \toprule
    module & $\hat H^0(G, \cdot)$ & $H^1(G,\cdot)$ & Conditions \\
    \midrule
    $\cyclic{p^i}$ & $\FF_p$ & $\FF_p$ & $i\geq 1$ \\
    $F_{i,\bar{w}}$ &$0$ & $0$ & $i\geq 2$ and $\bar{w}\in \FF_p^{\times}$ \\
    $\ZZ_p$ & $\FF_p$ & $0$ & \\
    $A$ & $0$ & $\FF_p$ & \\
    $\ZZ_p[G]$ & $0$ & $0$ & \\
    $\bigl\{\cyclic{p^i} \bigm \vert \ZZ_p\bigr\}$ & $\FF_p$ & $0$ & $i\geq 1$\\
    $\bigl\{\cyclic{p^i} \bigm \vert A \bigr\}$ & $0$ & $\FF_p$& $i\geq 1$\\
    $\bigl\{\cyclic{p^i} \bigm \vert \ZZ_p\oplus A \bigr\}$ & $0$ & $0$ & $i\geq 1$\\
  \end{tabular}\end{center}
  The decomposition is unique up to reordering the summands.
\end{thm}
\begin{proof}
  Let $M$ be a finitely generated $\ZZ_p[G]$-module.
  There is an exact sequence
  \begin{equation}\label{mt_seq}
    \xymatrix@1{0\ar[r] & M_t \ar[r] & M \ar[r] & M_f \ar[r] & 0 }
  \end{equation}
  where $M_f=M/M_t$ is a free, finitely generated $\ZZ_p$-module with an action by $G$.
  We assume that $M_t$ is a cyclic group, which implies by Lemma~\ref{cyclic_fin_lem} that $M_t$ is either $\cyclic{p^i}$ or $F_{i,\bar{w}}$ for some $i\geq 0$ and $\bar{w}\in\FF_p^{\times}$.
  Further by Proposition~\ref{cyclic_lattices_prop}, the lattice $M_f$ is a sum of $\ZZ_p$, $A$ and $\ZZ_p[G]$.

  We are now going to prove the following three statements:
  \begin{equation}\label{ext_eq}
    \Ext^1_G\bigl(M_f,F_{i,\bar{w}}\bigr)=0,\quad
    \Ext^1_G\bigl(\ZZ_p,\cyclic{p^i}\bigr) \cong \FF_p,\quad\text{ and }
    \Ext^1_G\bigl(A,\cyclic{p^i}\bigr) \cong \FF_p.
  \end{equation}
  It is clear that $\Ext^1_G(\ZZ_p[G],F_{i,\bar{w}})=0$ as $\ZZ_p[G]$ is free and hence projective.
  The norm map on $F_{i,\bar{w}}$ with $w= 1+p^{i-1}\bar{w}$ is the multiplication by
  \begin{equation*}
    1 + w + w^2+\cdots +w^{p-1} \equiv p + p^{i-1}\frac{p(p-1)}{2}\,\bar{w} \equiv p \pmod{p^i}
  \end{equation*}
  as $p$ is odd.
  Therefore $\hat H^0\bigl(G,F_{i,\bar{w}}\bigr) = F_{i,\bar{w}}^G/\No(F_{i,\bar{w}}) \cong pF_{i,\bar{w}}/pF_{i,\bar{w}}=0$.
  By the Herbrand quotient on finite modules $H^1\bigl(G,F_{i,\bar{w}}\bigr)=0$ and hence $\Ext^1_G\bigl(\ZZ_p, F_{i,\bar{w}}\bigr)$ vanishes, too.

  Consider the short exact sequence $\xymatrix@1{0\ar[r] &\ZZ_p\ar[r]^N &\ZZ_p[G] \ar[r] & A \ar[r] &0}$ which yields
  \begin{equation*}
    \xymatrix@1{0&\ar[l] \Ext^1_G\bigl(A,F_{i,\bar{w}}\bigr) & \ar[l] \Hom_G\bigl(\ZZ_p,F_{i,\bar{w}}\bigr) & \ar[l] \Hom_G\bigl(\ZZ_p[G], F_{i,\bar{w}}\bigr).}
  \end{equation*}
  By evaluation at $1$ the map on the right identifies with $\xymatrix@1{pF_{i,\bar{w}} &\ar[l]_{[p]} F_{i,\bar{w}}}$ and hence $\Ext^1_G\bigl(A,F_{i,\bar{w}}\bigr)=0$.
  This concludes the proof for the first statement in~\eqref{ext_eq}.

  Next, we see that $\Ext^1_G\bigl(\ZZ_p,\cyclic{p^i}\bigr) \cong H^1\bigl(G,\cyclic{p^i}\bigr) \cong \cyclic{p}$.
  Also with the same method as above, we find $\Ext^1_G\bigl(A,\cyclic{p^i}\bigr) \cong \cyclic{p}$ as it is the cokernel of the map $[p]$ on $\cyclic{p^i}$.
  This completes the proof of all statements in~\eqref{ext_eq}.

  We conclude that the only direct summands possibly appearing in $M$ are $\cyclic{p^i}$, $F_{i,\bar{w}}$, $\ZZ_p$, $A$, $\ZZ_p[G]$ or non-split exact sequences
  \[
    \xymatrix@1{0\ar[r]& \cyclic{p^i} \ar[r] & M \ar[r] & \ZZ_p^a \oplus A^b \ar[r] & 0}
  \]
  with $a,b\geq 0$.
  By the above such short exact sequences are parametrised by $\Ext^1_G\bigl(\ZZ_p^a\oplus A^b,\cyclic{p^i}\bigr)\cong \FF_p^{a+b}$.
  We are now going to show that there are only three distinct isomorphism classes of indecomposable $\ZZ_p[G]$-modules among these non-split exact sequences.

  First, we proceed to construct the extensions explicitly.
  For any $u\in\Hom\bigl(\ZZ_p^a,\FF_p\bigr)$, we set $M_u$ to be $\cyclic{p^i}\oplus \ZZ_p^a$ as a $\ZZ_p$-module, but with the action by a generator $\tau\in G$ defined by
  \[\tau\cdot(t,x) = (t+u(x)\,p^{i-1}, x)\text{ for }t\in\cyclic{p^i}\text{ and }x\in\ZZ_p^a.\]
  Then
  \begin{equation*}
    \xymatrix@1{ 0\ar[r] & \cyclic{p^i} \ar[r] & M_u \ar[r] & \ZZ_p^a \ar[r] & 0 }
  \end{equation*}
  is an non-split extension of $\ZZ_p[G]$-modules.

  The connecting homomorphism $\ZZ_p^a\to H^1\bigl(G,\cyclic{p^i}\bigr)\approx \cyclic{p}$ is equal to $u$, which shows that these are distinct extensions and they are non-split when $u\neq 0$.
  Assume $u\neq 0$, we find that $M_u^{G}$ consists of all $(t,x)$ with $x\in \ker (u)$.
  The norm map on $M_u$ sends $(t,x)$ to
  \begin{equation*}
    \No(t,x) = \sum_{j=0}^{p-1} \bigl(t + j\cdot u(x)\,p^{i-1},x\bigr) = \bigl(pt+p \tfrac{p-1}{2} u(x)p^{i-1},px\bigr)=\bigl(pt,px\bigr)
  \end{equation*}
  as $p$ is odd. Therefore $\hat{H}^0(G,M_u) \cong \FF_p^a$ and $H^1(G,M_u)=0$.

  Next, for any $\ZZ_p$-linear $v:\bigl(A/(\tau-1)A)^b\to\cyclic{p}$, we will build an extension
  \begin{equation*}
    \xymatrix@1{ 0\ar[r] & \cyclic{p^i} \ar[r] & M'_v \ar[r] & A^b\ar[r] & 0}
  \end{equation*}
  as follows.
  We identify $A$ with $\ZZ_p[\zeta]$ where $\zeta^p=1$ and the action by $\tau$ is by multiplication with $\zeta$.
  We define $M'_v$ as the $\ZZ_p$-module $\cyclic{p^i}\oplus \ZZ_p[\zeta]^b$ together with the action by $\tau$ given by $\tau(t,y)=(t + f_v(y),\zeta\cdot y)$ where $f_v$ is a $\ZZ_p$-linear map $A^b\to \cyclic{p^i}$ such that $f_v\bigl( \tfrac{p}{1-\zeta} y\bigr)$ reduces to $v(y)$ modulo $p$ for all $y$.
  This $f_v$ exists as the condition only imposes the values on a $b$ dimensional $\FF_p$-subspace of the $b(p-1)$ dimensional space $(A/pA)^b$.
  We find for $0\leq j < p$
  \begin{equation*}
    \tau^j\,(t,y) = \bigl(t + f_v(y) + f_v(\zeta\,y) + \cdots+f_v(\zeta^{j-1} y), \zeta^j y\bigr) = \bigl( t + f_v ( \tfrac{1-\zeta^j}{1-\zeta} y),\zeta^j y\bigr).
  \end{equation*}
  and therefore
  \begin{align*}
    N(t,y) &= \biggl( pt + \sum_{j=0}^{p-1} f_v\Bigl(  \tfrac{1-\zeta^j}{1-\zeta} y\Bigr),\ \sum_{j=0}^{p-1} \zeta^j y\biggr)
    \\ &= \biggl( pt + f_v\Bigl( \sum_{j=0}^{p-1} (1-\zeta^j) \tfrac{y}{1-\zeta}\Bigr),\ 0\biggr)
   = \Bigl( pt + f_v(\tfrac{p}{1-\zeta}y),\ 0\Bigr)
  \end{align*}
  The connecting homomorphism $\hat H^{-1}(G,A^b)\to \hat{H}^0\bigl(G,\cyclic{p^i}\bigr)$ identifies with $v$.
  When $v\neq 0$, then $\hat H^0(G,M_v)=0$ and $H^1(G,M_v)\cong\FF_p^b$.

  Finally, the extension $M_{u,v}$ is defined as the $\ZZ_p$-module $\cyclic{p^i}\oplus \ZZ_p^a\oplus A^b$ with $\tau$ acting on $(t,x,y)$ by $(t+u(x)p^{i-1}+f_v(y),x,\zeta\cdot y)$ with $u$, $v$ and $f_v$ as above.
  It is not hard to calculate that $\dim_{\FF_p} \hat H^0\bigl(G,M_{u,v}\bigr) = a-1$ and $\dim_{\FF_p} H^1\bigl(G,M_{u,v}\bigr)=b-1$.

  With the above, we have explicitly constructed all classes of extensions in
  \[\Ext^1_G\bigl(\ZZ_p^a\oplus A^b,\cyclic{p^i}\bigr) \cong\Hom\bigl(\ZZ_p^a,\FF_p\bigr)\oplus \Hom_{\ZZ_p[\zeta]}\bigl(\ZZ_p[\zeta], \FF_p\bigr) \cong \FF_p^a\oplus \FF_p^b.\]

  The group $\Aut_G\bigl(\ZZ_p^a\oplus A^b\bigr)$ acts from the left on this extension group and this action does not change the isomorphism class of $M_{u,v}$ as a $\ZZ_p[G]$-module.
  The group acting is isomorphic to $\GL_a(\ZZ_p)\times \GL_b(\ZZ_p[\zeta])$ and, for $\alpha\in\GL_a(\ZZ_p)$ and $\beta\in\GL_b(\ZZ_p[\zeta])$, the action on $(u,v) $ gives $\bigl( u\circ \alpha^{-1}, v\circ \beta^{-1}\bigr)$.
  It follows that there are four orbits on $\FF_p^a\oplus \FF_p^b$ corresponding to $u$ and $v$ being zero or non-zero.
  We set $\bigl\{\cyclic{p^i}\bigm\vert \ZZ_p\bigr\} := M_u$, $\bigl\{\cyclic{p^i}\bigm\vert A\bigr\} := M_v$ and $\bigl\{\cyclic{p^i}\bigm\vert \ZZ_p\oplus A\bigr\} := M_{u,v}$ for any choice of non-zero $u$ and $v$ with $a=b=1$.
  These three represent distinct isomorphism classes of indecomposable $\ZZ_p[G]$-modules.
  For a general $u$ or $v$ with $a>1$ or $b>1$, we can show that the $\ZZ_p[G]$-module $M_{u,v}$ is isomorphic to a direct sum of one of the above and a sum $\ZZ_p^{a-1}\oplus A^{b-1}$ in the following way illustrated in the case $a>1$ and $b=0$:

  Suppose $0\neq u\in\Hom\bigl(\ZZ_p^a,\FF_p\bigr)$ and $a>1$.
  Then there is a linear transformation $\alpha\in\GL_a(\ZZ_p)$ such that $u'=u\circ \alpha^{-1}$ maps the first basis vector $e_1 \in \ZZ_p^a$ to $1$ and all the other basis vectors $e_i$ with $i\neq 0$ to~$0$.
  The map $M_u \to M_{u'}$ sending $(t,x)$ to $(t,\alpha(x))$ is the $\ZZ_p[G]$-isomorphism mentioned above.
  Since the action of $G$ for $M_{u'}$ only involves the first basis vector we see that $M_{u'} \cong \bigl\{\cyclic{p^i}\bigm\vert \ZZ_p\bigr\} \oplus \ZZ_p^{a-1}$.
  This argument also shows directly that the choice of $u$ in the definition of $\bigl\{\cyclic{p^i}\bigm\vert \ZZ_p\bigr\}$ does not change the isomorphism class of the $\ZZ_p[G]$-module.
  The situation for $M_v$ and $M_{u,v}$ is analogous.

  The last statement, that the direct sum is unique, is a consequence of the Krull-Schmidt-Azumaya Theorem.
\end{proof}

\begin{lem}\label{str_lem}
  There are non-trivial extensions $\xymatrix@1{0\ar[r]&\ZZ_p[G] \ar[r]& M \ar[r] & F \ar[r] & 0}$ where $F$ is one of the following finite $G$-modules: $F=\cyclic{p}$, $F=\FF_p[ G]$ or $F\cong F_{2,\bar{w}}$. Moreover
  \begin{enumerate}
    \item If $F=\cyclic{p}$ then $M\cong \ZZ_p\oplus A$ as a $\ZZ_p[G]$-module.
    \item If $F= F_{2,\bar{w}}$ or if $F=\FF_p[ G]$ and $M_{t}$ is cyclic, then either $M\cong \ZZ_p[G]$ or $M\cong\bigl\{\cyclic{p}\bigm\vert\ZZ_p\oplus A\bigr\}$.
  \end{enumerate}
\end{lem}
\begin{proof}
  For the first statement, we need to compute $\Ext^1_G\bigl(F,\ZZ_p[G]\bigr)$.
  Let $p^k$ be the exponent of $F$ as an abelian group and consider the multiplication by $p^k$ on $\ZZ_p[G]$.
  From
  \begin{equation*}
    \xymatrix@1{\Hom_G\bigl(F,\ZZ_p[G]\bigr)=0\ar[r] & \Hom_G\bigl(F,\cyclic{p^k}[ G]\bigr)\ar[r] & \Ext^1_G\bigl(F,\ZZ_p[G]\bigr) \ar^{[p^k]=0}[r] & \Ext^1_G\bigl(F,\ZZ_p[G]\bigr) }
  \end{equation*}
  we see that $\Ext^1_G\bigl(F,\ZZ_p[G]\bigr) \cong \Hom_G\bigl(F, \cyclic{p^k} [G]\bigr)$.
  First, the map from $F_{2,\bar{w}}$ to $\cyclic{p^2}[G]$ sending $1$ to $\sum_{i=0}^{p-1} w^{-i}\tau^i$ with $w=1+p\bar{w}$ is a $G$-equivariant map; therefore $\Ext^1_G\bigl(F_{2,\bar{w}},\ZZ_p[G]\bigr) \neq 0$.
  Since $\Hom_G\bigl(\cyclic{p}, \FF_p[G]\bigr) \cong \cyclic{p}$ and $\Hom_G\bigl(\FF_p[G],\FF_p[G]\bigr) \neq 0$ we also have the non-trivial extensions in the other cases.

  Assume now $M$ is such a non-trivial extension.
  Since $\ZZ_p[G]$ is cohomologically trivial, we have $H^i(G,M)\cong H^i(G,F)$.
  Note also that $M_{t}$ injects into $F$, but it cannot surject otherwise the extension would be split.

  For $F=\cyclic{p}$ the non-trivial extensions must then be torsion-free.
  Since $H^1\bigl(G,\cyclic{p}\bigr)\cong\cyclic{p}$, the table in Theorem~\ref{cyclic_class_thm} tells us that $M$ must be isomorphic to $\ZZ_p \oplus A$.

  If $F$ is $F_{2,\bar{w}}$, then $M_{t}$ is either trivial or $\cyclic{p}$.
  Since $H^1\bigl(G,F_{2,\bar{w}}\bigr)=0$, this only leaves $\ZZ_p[G]$ or $\bigl\{\cyclic{p}\bigm\vert \ZZ_p\oplus A\bigr\}$; in the first case $1\in\ZZ_p[G]$ is sent to $\tau-\tilde{w}$, where $\tilde{w}\in\ZZ_p^{\times}$ is a lift of $w$, in the second case it is sent to $(0,p,1)\in \bigl\{\cyclic{p}\bigm\vert \ZZ_p\oplus A\bigr\}$ defined as $M_{u,v}$ in the proof of Theorem~\ref{cyclic_class_thm}.

  If $F=\FF_p[G]$ and $M_{t}$ is cyclic then $M_{t}$ has at most $p$ elements.
  We reach the same conclusion as above given that $H^1\bigl(G,\FF_p[G]\bigr)$ also vanishes.
  Also here it is possible to write down an explicit extension.
\end{proof}

\subsection{Metacyclic groups}\label{metacyc_subsec}

Despite being interested mainly in the dihedral case, we present the results in a slightly more general setting.
Let $m>1$ be a divisor of $p-1$.
Let $\xi$ be a choice of a primitive $m$-th root of unity in $\ZZ_p$.
Let $r\in \ZZ$ such that $r\equiv \xi\pmod{p}$ so that $r^m\equiv 1 \pmod{p}$.
In this section we treat the case of the metacyclic group
\begin{equation}\label{metacyc_eq}
 G = \bigl\langle \tau,\sigma \bigm\vert \tau^p = \sigma^m = 1 \text{ and }\sigma\tau = \tau^r\sigma\bigr\rangle
\end{equation}
of order $pm$.
We write $N$ for the normal subgroup generated by $\tau$ and $H$ for the subgroup generated by $\sigma$; so $G = N \rtimes H$.
This includes the case $G$ is the dihedral group $D_p$ when $m=2$ and $r=-1$.

For $i\in \cyclic{m}$, we define the $G$-lattice $\ZZ_p\{i\}$ which, as a group, is just $\ZZ_p$, the action by $\tau$ is trivial and by $\sigma$ is the multiplication with $\xi^i$.
We have $\ZZ_p\bigl[G/N\bigr] \cong \bigoplus_{i=0}^{m-1} \ZZ_p\{i\}$ as the index of $N$ in $G$ is coprime to $p$.

Next, let $B=\ZZ_p\bigl[G/H\bigr]$ and set $B\{i\} = B \otimes_{\ZZ_p} \ZZ_p\{i\}$. Then $\ZZ_p[G] \cong \bigoplus_{i=0}^{m-1} B\{i\}$.

Finally, we define $A$ as $(\tau-1)B$, which is the kernel of the augmentation map $B\to \ZZ_p$.
If we equip the group $\ZZ_p[\zeta]$ where $\zeta$ is a primitive $p$-th root of unity with the group action defined by letting $\tau$ act as multiplication by $\zeta$ and $\sigma(\zeta^j) = \zeta^{rj}$ for all $0\leq j <p$, then $A$ is isomorphic to $(\zeta-1)\ZZ_p[\zeta]$.
For any $i$, we set $A\{i\} = A \otimes_{\ZZ_p} \ZZ_p\{i\}$.
In particular $A\{m-1\} \cong \ZZ_p[\zeta]$.
It follows that $\Ind_{N}^G (A) \cong \bigoplus_{i=0}^{m-1} A\{i\}$.

We have a non-split exact sequences
\[
 \xymatrix@R=1ex{
 0\ar[r]& A\{i\}\ar[r] &B\{i\}\ar[r] & \ZZ_p\{i\}\ar[r] & 0 \\
 0\ar[r]& \ZZ_p\{i\}\ar[r] &B\{i\}\ar[r] & A\{i-1\}\ar[r] & 0
 }
\]
for all $0\leq i <m$.

The complete list of simple $\QQ_p[G]$-modules is given by $\ZZ_p\{i\}\otimes\QQ_p$ and $A\{i\}\otimes \QQ_p$; this follows from Proposition~\ref{pu_prop} below.
It was first proven by Pu in~\cite{pu}.
See~\cite{lee} for the case of dihedral groups.

\begin{prop}\label{pu_prop}
 The lattices $\ZZ_p\{i\}$, $A\{i\}$ and $B\{i\}$ for $0\leq i <m$ represent all isomorphism classes of indecomposable finitely generated $\ZZ_p[G]$-lattices.
\end{prop}
\begin{proof}
 Since $[G:N]$ is coprime to $p$, the proof of Proposition~33.4 in~\cite{cr1} shows that all indecomposable modules are summands of $\Ind_N^G(X)$ as $X$ runs through all indecomposable $N$-lattices.
 The proposition follows now from Proposition~\ref{cyclic_lattices_prop} and $\Ind_N^G(\ZZ_p)=\bigoplus_{i=0}^{m-1} \ZZ_p\{i\}$, $\Ind_N^G(A) = \bigoplus_{i=0}^{m-1} A\{i\}$ and $\Ind_N^G(\ZZ_p[N])= \bigoplus_{i=0}^{m-1} B\{i\}$.
\end{proof}

As all $B\{i\}$ are direct summands of the free $\ZZ_p[G]$, they are projective $G$-lattices.
Therefore they are cohomologically trivial.
Write $\FF_p\{i\}$ for the $G$-module $\ZZ_p\{i\}/p\ZZ_p\{i\}$.
Then we have that $\hat H^0\bigl(N,\ZZ_p\{i\}\bigr)\cong\FF_p\{i\}$ and $ H^1\bigl(N,A\{i\}\bigr)\cong \FF_p\{i\}$ as $G/N$-modules.

Unlike for the cyclic group, the decomposition of a $\ZZ_p[G]$-lattice $M$ into indecomposable lattices cannot be determined by knowing the $\QQ_p[G]$-module $M\otimes \QQ_p$ and the cohomology groups $\hat H^i\bigl(N,M)$ only.
In~\cite{torzewski}, Torzewski considered an additional invariant, which we are going to introduce next.

For a general finite group $\mathcal{G}$, we say that a $\ZZ_p[\mathcal G]$-lattice $M$ is rationally self-dual if $\Hom_{\QQ_p}\bigl(M\otimes\QQ_p,\QQ_p) \cong M\otimes\QQ_p$ as $\QQ_p[\mathcal G]$-modules; equivalently there is a non-degenerate $\mathcal G$-equivariant symmetric bilinear pairing $\beta\colon M \times M \to \QQ_p$.
If $M\cong \tilde M \otimes_{\ZZ}\ZZ_p$ for some $\ZZ[\mathcal G]$-lattice $\tilde M$ then it is automatically rationally self-dual.
All $\ZZ_p[\mathcal G]$-lattices are rationally self-dual if $\mathcal G$ is the cyclic group of order $p$ or the dihedral group of order $2p$.
Instead for our group $G$, any rationally self-dual lattice is a (not necessarily unique) direct sum of
\begin{gather*}
  \ZZ_p,\ A,\ B,\ A\{i\},\ \ZZ_p\{i\}\oplus \ZZ_p\{-i\},\ B\{i\}\oplus \ZZ_p\{-i\},\ B\{i\}\oplus B\{-i\} \text{ for $0<i<\tfrac{m}{2}$ }\\
  \text{ as well as } B\{\tfrac{m}{2}\},\ \ZZ_p\{\tfrac{m}{2}\}\text{ when $m$ is even.}
\end{gather*}

A Brauer relation for $G$ is a formal sum $\Theta = \sum_{H} a_H\, H$ of subgroups $H\leq G$ with coefficients $a_H\in\ZZ$ such that $\bigoplus_H \QQ\bigl[G/H\bigr]^{a_H}$ is zero as a virtual representation of $G$.
For our group $G$, all Brauer relations can be obtained from Artin's induction theorem, see~Theorem~2.10 in~\cite{torzewski}:
For each divisor $1<d\mid m$, let
\[
  \Theta_d = 1 - d\cdot H_d - N + d\cdot G_d
\]
where $H_d$ is the subgroup of $H$ of order $d$ and $G_d = N\rtimes H_d$.

To every Brauer relation $\Theta=\sum_H a_H\, H$ and every rationally self-dual $\ZZ_p[G]$-lattice $M$ one associates a Dokchitser regulator constant $\mathcal{C}_{\Theta}(M)\in \QQ_p^{\times}/\square$ where $\square =\bigl\{z^2\bigm\vert z\in\ZZ_p^{\times}\bigr\}$ by
\begin{equation}\label{reg_eq}
  \mathcal{C}_{\Theta}(M) = \prod_{H} \det\Bigl( \frac{1}{\vert H \vert} \beta\Bigm \vert M^H \Bigr)^{a_H} \cdot \square
\end{equation}
where $\beta$ is a $G$-equivariant pairing on $M$.
See~\cite{DD09} for the basic properties, including the fact that the definition is independent of the choice of $\beta$.
The special case of dihedral groups was already worked out by Bartel in Theorem~4.4 in~\cite{bartel}.

We are interested in the integer $s_d(M)$ defined to be the $p$-adic valuation of $\mathcal{C}_{\Theta_d}(M)$.
As proved in~\cite[Theorem 1.1]{torzewski}, the kernel of the map $s = \oplus s_d$ from the $\QQ$-vector space with  basis $\ZZ_p[G/U]$ as $U$ runs through all subgroups $U\leq G$ to $\QQ$ is equal to the subspace generated by all cyclic $U$.

\begin{prop}
  Let $1<d$ be a divisor of $m$ and $0<i<m$.
  We have $s_d\bigl(\ZZ_p\bigr)=1-d$ and $s_d\bigl(\ZZ_p\{i\}\oplus\ZZ_p\{-i\}\bigr) = 2$.
  Also $s_d(A)=1-d$ and $s_d\bigl(A\{i\}\bigr)=s_d\bigl(B\{i\}\oplus \ZZ_p\{-i\}\bigr)=2i+1-d$ and $s_d\bigl(B\bigr) = s_d\bigl(B\{i\}\oplus B\{-i\}\bigr) = 0$.
  Further, if $m$ is even, $s_d\bigl(\ZZ_p\{\tfrac{m}{2}\}\bigr) = 1$ and $s_d\bigl(B\{\tfrac{m}{2}\}\bigr) = 0$.
\end{prop}
Note that $s_d$ cannot be extended to an additive function on all lattices.
The result in the proposition can be deduced from the explicit and more general calculation by Torzewski in~\cite[Proof of Theorem~4.1]{torzewski} where he found that $s_{d}\bigl( \ZZ_p[G/G_e] \bigr) = (1-g)\cdot \tfrac{m\,g}{d\,e}$ for any $d$ and $e$ dividing $m$ and $g=\gcd(d,e)$.
We proceed here to calculate it directly on all the minimal rationally self-dual lattices.

\begin{proof}
  By Proposition~2.45.(3) in~\cite{DD09}, the regulator constants satisfy $\mathcal{C}_{\Theta_d}(M) = \mathcal{C}_{\Theta_d}\bigl( M\vert_{G_d} \bigr)$, which implies that we may restrict to the case when $d=m$.
  We will write $\mathcal{C}$ for $\mathcal{C}_{\Theta_m}$.
  To show that $\mathcal{C} \bigl(\ZZ_p\bigr)=p^{1-m}$ is a direct calculation on the definition.
  The same goes for $\mathcal{C}\bigl(\ZZ_p\{i\}\oplus\ZZ_p\{-i\}\bigr) = p^2$ using the bilinear form $\bigr(\begin{smallmatrix} 0 & 1\\ 1& 0 \end{smallmatrix}\bigr)$ for all $0<i<m$.
  From the additivity of the regulator constant, see Corollary~2.18 in~\cite{DD09}, we may conclude that $\mathcal{C}\bigl(\ZZ_p\{\tfrac{m}{2}\}\bigr) = \pm p$ for even $m$.

  Let $y = \sum_{j=0}^{m-1} \xi^{-j} \tau^{r^j}H\in A$.
  It satisfies $\sigma(y) = \xi\cdot y$.
  Using that $\xi\equiv r \pmod{p}$, we find
  \[
   y  \equiv \Bigr(\sum_{j=0}^{m-1} \xi^{-j}\Bigl) + \Bigl( \sum_{j=0}^{m-1} \xi^{-j} r^j \Bigr)\cdot (\tau - 1 )H
      \equiv 0 + m \cdot (\tau-1) \not\equiv 0 \pmod{ (\tau-1)^2B }.
  \]
  Consider the map $A\otimes \ZZ_p\{i+1\}\to A\otimes \ZZ_p\{i\}$ sending $a\otimes z$ to $ay\otimes z$.
  Its image is $(\tau-1)A\otimes \ZZ_p\{i\}$, which is of index $p$.
  Therefore we have the short exact sequence of $\ZZ_p[G]$-modules
  \begin{equation}\label{as_seq}
    \xymatrix@1{ 0 \ar[r] &A\{i+1\} \ar[r] & A\{i\} \ar[r] & \FF_p\{i+1\}\ar[r] & 0. }
  \end{equation}
  We conclude that $A\{i+1\}$ has index $p$ in $A\{i\}$ and that $A\{i+1\}^H = A\{i\}^H$ for all $0\leq i<m-1$.
  Using $A\{i\}^N=A\{i\}^G=0$, for all $0\leq i<m-1$, we have
  \[
    \mathcal{C}\bigl(A\{i+1\}\bigr) = \frac{ \det\bigl( \beta \bigm \vert A\{i+1\} \bigr) }{\det \bigl( \tfrac{1}{m} \beta \bigm\vert A\{i+1\} ^H\bigr)^m }
    = \frac{ p^2 \det\bigl( \beta \bigm \vert A\{i\} \bigr) }{\det \bigl( \tfrac{1}{m} \beta \bigm\vert A\{i\}^H \bigr)^m } = p^2\cdot \mathcal{C}\bigl(A\{i\}\bigr).
  \]
  Thus $\mathcal{C} \bigl(A\{i\}\bigr) = p^{2i}\cdot\mathcal{C}( A )$ for all $0\leq i\leq m-1$.
  Since $\bigoplus_{i=0}^{m-1} A\{i\} \cong \Ind_N^G(A)$, we find
  \[
    \mathcal{C}\bigl( \Ind_N^G(A) \bigr) = \prod_{i=1}^m \mathcal{C}\bigl(A\{i\}\bigr) = \mathcal{C}\bigl(A\{1\}\bigr)^m \cdot p^{m(m-1)}.
  \]
  But since $\mathcal{C}_{\Theta_m}\bigl( \Ind_N^G(A) \bigr) = \mathcal{C}_{\res_N(\Theta_m)}\bigl(A\bigr)$ by Proposition~2.45 in~\cite{DD09} and the restriction of $\Theta_m$ to $N$ is trivial, we get $\mathcal{C}(A)^m = p^{-m(m-1)}$.
  Hence $s_m(A) = 1-m$ and by the above recursion formula $s_m\bigl(A\{i\}\bigr) = 2\,i +1-m$ for all $0\leq i \leq m-1$.
  In particular $s_m\bigl(\ZZ_p[\zeta]\bigr) = m-1$.

  Just as before $\mathcal{C}\bigl(B\bigr) = \mathcal{C}\bigl(\ZZ_p[G/H]\bigr) = 1$ because the restriction of $\Theta_m$ to $H$ is trivial; that is Lemma~2.46 in~\cite{DD09}.

  Consider the map
  \begin{align*}
    \Phi\colon \ZZ_p\bigl[G/H\bigr] &\to \ZZ_p\oplus \ZZ_p[\zeta] \\
    \sum_{i=0}^{p-1} a_i \,\tau^i H  &\mapsto \Bigl( \sum_{i=0}^{p-1} a_i, \ \sum_{i=0}^{p-1} a_i \zeta^i \Bigr).
  \end{align*}
  It is an injective $G$-equivariant map whose image is $\bigl\{(x,a) \in \ZZ_p\oplus \ZZ_p[\zeta] \bigm\vert x\equiv a \pmod{\zeta-1}\bigr\}$ of index~$p$.
  We deduce the exact sequence
  \[
    \xymatrix@1{ 0\ar[r] & B\{i\}\oplus \ZZ_p\{-i\} \ar[r] & A\{i-1\} \oplus \ZZ_p\{i\}\oplus \ZZ_p\{-i\} \ar[r] & \FF_p\{i\}\ar[r] & 0 }
  \]
  for all $i$.
  When $0<i<m$, then the $G$ and $H$-invariant parts of the first two terms are equal, while the $N$-invariant parts are of index $p$.
  The same reasoning as above concludes now that
  \[
\mathcal{C}\bigl(B\{i\}\oplus \ZZ_p\{-i\} \bigr) = \mathcal{C}\bigl(A\{i-1\} \oplus \ZZ_p\{i\}\oplus \ZZ_p\{-i\}\bigr) =
\mathcal{C}\bigl(A\{i-1\}\bigr) \cdot \mathcal{C}\bigr(\ZZ_p\{i\}\oplus \ZZ_p\{-i\}\bigr)
  \]
  and hence $s_m\bigl(B\{i\}\oplus \ZZ_p\{-i\} \bigr) = 2(i-1)+1-m+2=2\,i+1-m$ for $1\leq i \leq m-1$.

  Finally
  \begin{align*}
    \mathcal{C}\bigl(B\{i\}\oplus B\{-i\}\bigr)\cdot \mathcal{C}\bigl( \ZZ_p\{i\}\oplus\ZZ_p\{-i\}\bigr) &=
\mathcal{C}\bigl(B\{i\}\oplus B\{-i\}\oplus \ZZ_p\{i\}\oplus\ZZ_p\{-i\}\bigr) \\
&= \mathcal{C}\bigl(B\{i\}\oplus \ZZ_p\{-i\}\bigr)\cdot \mathcal{C}\bigl( \ZZ_p\{i\}\oplus B\{-i\}\bigr)
  \end{align*}
  shows that, for $0<i<m$, we have $s_d\bigl(B\{i\}\oplus B\{-i\}\bigr) = 2\,i+1-m + 2(m-i)+1-m - 2 = 0$.
\end{proof}

In the special case $m=2$, when the group is dihedral, all $G$-lattices are rationally self-dual.
\begin{thm}\label{torzewski_thm}
  Let $p>2$ be a prime and let $G=D_p$.
  Then a $\ZZ_p[G]$-lattice~$M$ is determined up to isomorphism by the knowledge of
  \begin{itemize}[nosep]
    \item $M\otimes \QQ_p$ as a $\QQ_p[G]$-module,
    \item $H^1(N,M)$ as a $\FF_p[G/N]$-module, and
    \item the regulator constant $s(M) = \ord_p\bigl( \mathcal{C}_{\Theta_2}(M)\bigr)$.
  \end{itemize}
\end{thm}
This is explained in~\cite[Section 7.1]{torzewski}, but can also be read off the Table~\ref{dihedral_table} below.
This theorem does not extend to the more general metacyclic groups even if one restricts to rationally self-dual $G$-lattices.

\subsection{Saturation index}\label{saturation_subsec}

Given a finite group $G$ and a $\ZZ_p[G]$-module $M$, we define the saturation index $\iota(M)$ to be the quotient of $M$ by the subgroup generated by all $M^H$ where $H$ runs through all non-trivial cyclic subgroups of $G$.
Alternatively, it is the quotient of $M$ by the sub-$\ZZ_p[G]$-module generated by all $M^H$ where $H$ runs through a set of representatives of all conjugacy classes of non-trivial cyclic subgroups of $G$.
\[
  \iota(M) = \frac{M}{\sum_{\text{cyclic }H\leq G} M^H}. 
\]
This index and its generalisations appear dominantly in the work of Bartel and de Smit~\cite{bartel_desmit, bartel}.
By Artin's induction theorem $\iota(M)$ is a finite $\ZZ_p[G]$-module for all non-cyclic groups $G$.

We call $\iota(M)$ the saturation index because of the following observation in the case of $M$ being the $p$-adic completion of $E(K)$ for some elliptic curve. If we need to determine $E(K)$ explicitly, then we would start by a search for points in $E(K^H)$ for proper subgroups $H$ as this is quicker than searching in $E(K)$ directly.
These points then generate a submodule of $E(K)$ of finite index if $G$ is non-cyclic and there is an effective algorithm~\cite{prickett}, called a $p$-saturation, to calculate the full $M$.
This algorithm effectively calculates $\iota(M)$, though usually only its size is of interest.
Unfortunately, at this stage, we have no means to relate the invariant $\iota(M)$ directly to arithmetic information of $E$ that is easier to calculate than $M$ itself.
We may use it to determine $M$ once we have $E(K^H)$ for all non-trivial cyclic $H$, without having to calculate the matrices representing the action of $G$ on the generators of $E(K)$.

Note also that the functor $\iota$ is additive $\iota(M\oplus M')=\iota(M)\oplus\iota(M')$ because $\sum_H M^H \oplus M'^H = \bigl(\sum_H M^H\bigr) \oplus \bigl(\sum_H M'^H\bigr)$, but it does not behave well in short exact sequences.

\begin{prop}\label{saturation_prop}
  Let $G$ be a metacyclic group of order $pm$ as in~\eqref{metacyc_eq}.
  Then $\iota(M)$ is trivial for $M$ isomorphic to $A\{-1\}\cong\ZZ_p[\zeta]$, $B$ or $\ZZ_p\{i\}$ for any $i$, while $\iota \bigl( A\{i-1\}\bigr) \cong \iota \bigl( B\{i \} \bigr)$ is a $\FF_p$-vector space of dimension $m-i$  for any $0<i<m$.
\end{prop}

The structure of $\bigl( A\{i-1\}\bigr)$ as a $\ZZ_p[G]$-module can be described by stating that it is the unique indecomposable $\FF_p[N]$-module of that dimension and as a $\FF_p[H]$-module it is isomorphic to $\bigoplus_{k=1}^{m-i} \FF_p\{k\}$.
This extends the calculation by Bartel in~\cite{bartel} to more general metacyclic groups.

\begin{proof}
If $M=\ZZ_p\{i\}$ for any $i$ then $M^N = M$ and therefore $\iota(M)=0$.
Also $\iota(B)=0$ because the element $1\,H\in \ZZ_p\bigl[G/H\bigr]$ generates $B$ and it is fixed by the action of the cyclic group $H$.
Similarly, $\iota(\ZZ_p[\zeta])=0$ as $\ZZ_p[\zeta]^H$ contains $1$.

Let now $0<i<m$ and $M=A\{i-1\}\cong (\zeta-1)^i\,\ZZ_p[\zeta]$.
Since $M^N=0$, we only need to find $M^H$.
From sequence~\eqref{as_seq} we find that $A\{i+1\}^H$ is isomorphic to $A\{i\}^H$.
Hence by induction $\bigl((\zeta-1)^i\,\ZZ_p[\zeta]\bigr)^H = (\zeta-1)^{m} \ZZ_p[\zeta]\cap \ZZ_p[\zeta]^H$.
Therefore $\iota\bigl((\zeta-1)^i\ZZ_p[\zeta]\bigr)$ is the quotient of $(\zeta-1)^i$ by $(\zeta-1)^m$ as ideals in $\ZZ_p[\zeta]$.
Hence $\iota\bigl(A\{i-1\}\bigr)$ is isomorphic as a group to $\ZZ_p[\zeta]/(\zeta-1)^{m-i}$, which is isomorphic to $\FF_p^{m-i}$ as $m-i<p-1$

If we identify $B\{i\}$ with the subset of $(a,z)\in (\zeta-1)^i\ZZ_p[\zeta]\oplus\ZZ_p\{i\}$ such that $a\equiv (\zeta-1)^i\,z \pmod{(\zeta-1)^{i+1} }$, then we find $\bigl(B\{i\}\bigr)^N = \bigl\{ (0,z) \bigm\vert z\in p\ZZ_p\bigr\}$ and $\bigl(B\{i\}\bigr)^H = \bigl\{ (a,0) \bigm\vert a \in (\zeta-1)^{m}\ZZ_p[\zeta]\cap\ZZ_p[\zeta]^H\bigr\}$.
It follows that $\iota\bigl(B\{i\}\bigr)\cong \iota\bigl((\zeta-1)^i\ZZ_p[\zeta]\bigr)=\iota\bigl(A\{i-1\}\bigr)$.
\end{proof}

The following can be read out directly from Table~\ref{dihedral_table} below.

\begin{prop}
  The last entry in the list in Theorem~\ref{torzewski_thm} can be replaced by
  \begin{itemize}[nosep]
    \item $\dim_{\FF_p} \iota(M)$.
  \end{itemize}
\end{prop}

\subsection{Summary}

In Table~\ref{dihedral_table} we summarise the information gathered about $\ZZ_p[G]$-lattices, in case $G = D_p$.
We write $\twist{\ZZ}_p$ for $\ZZ_p\{1\}$, as well as $\twist{A} = A\{1\}$, and $\twist{B} = B\{1\}$, and $\twist{\FF}_p = \FF_p\{1\}$.
\begin{table}[ht]
  \centering
    \caption{Invariants of lattices for the dihedral group $D_p$}\label{dihedral_table}
    \renewcommand{\arraystretch}{1.1}
    \begin{tabular}{r*{7}{c}}
       \toprule
  $M$ && $\ZZ_p$ & $\twist{\ZZ}_p$ & $A$ & $\twist{A}$ & $B$ & $\twist{B}$ \\
       \midrule
  $\rk M$   && $1$ & $1$ & $p-1$ & $p-1$ & $p$ & $p$ \\
  $\rk M^H$ && $1$ & $0$ & $\tfrac{p-1}{2}$ & $\tfrac{p-1}{2}$ & $\tfrac{p+1}{2}$ & $\tfrac{p-1}{2}$ \\
  $\rk M^N$ && $1$ & $1$ & $0$ & $0$ & $1$ & $1$ \\
  $\rk M^G$ && $1$ & $0$ & $0$ & $0$ & $1$ & $0$ \\
  $\hat H^0(N,M)$ && $\FF_p$ & $\twist{\FF}_p$ & $0$ & $0$ & $0$ & $0$\\
  $     H^1(N,M)$ && $0$ & $0$ & $\FF_p$ & $\twist{\FF}_p$ & $0$ & $0$  \\
  $s(M)$ && $-1$ & $1$ & $-1$ & $1$ & $0$  & $0$ \\
  $\dim_{\FF_p} \iota(M)$ && $0$ & $0$ & $1$ & $0$ & $0$ & $1$\\
       \bottomrule
    \end{tabular}
\end{table}

\section{The group of local points in a cyclic extension}\label{local_mw_sec}

We now turn our attention to the $\ZZ_p[G]$-structure of the group of points on an elliptic curve in a local extension.
The method in this section could be applied to an arbitrary extension whose group is one we understand the $\ZZ_p[G]$-modules that could arise.
However we will concentrate on the simplest extensions and it turns out that the answer is already quite involved.

\begin{thm}\label{local_thm}
 Let $p>2$ be a prime and let $K/\QQ_p$ be the unramified extension of degree $p$.
 Let $E/\QQ_p$ be an elliptic curve.
 Suppose that, if $p=3$, the curve does not have additive reduction of type IV or IV${}^*$.
 Unless the reduction is split multiplicative and the Tamagawa number $c_v$ is divisible by $p$, we are in one of the following three cases:
 \begin{center}\begin{tabular}{rccc}
  \toprule
   $\vert E(\QQ_p)[p]\vert $ & $1$ & $p$ & $p$ \\
  $\vert E(K)[p]\vert $ & $1$ & $p$ & $p^2$ \\
  $E(K)\hatox \ZZ_p$ & $\ZZ_p[G]$ & $\bigl\{\ZZ/p\ZZ\bigm\vert \ZZ_p\oplus A\bigr\}$ & $F_{2,\bar{w}}\oplus \ZZ_p[G]$ \\
  \bottomrule
  \end{tabular}\end{center}
 If the reduction is split multiplicative and $p\mid c_v$, then set $j=\ord_p(c_v)$ and $\vert E(\QQ_p)[p^\infty]\vert = p^i$.
 Then we are in one of the following cases:
 \begin{center}\begin{tabular}{rccc}
  \toprule
   & $i=0$ & $ j=i>0$ & $j>i>0$ \\
   $E(K)\hatox \ZZ_p$ & $\ZZ_p\oplus A$ & $\cyclic{p^i}\oplus \ZZ_p[G]$ & $\bigl\{\cyclic{p^i}\bigm\vert A\bigr\}\oplus \ZZ_p$ \\
  \bottomrule
 \end{tabular}\end{center}
\end{thm}

Note that the modules $\bigl\{\cyclic{p^i}\bigm\vert \ZZ_p\oplus A\bigr\}$ with $i>1$ and the modules $\bigl\{\cyclic{p^i}\bigm\vert \ZZ_p\bigr\}\oplus A$ for $i\geq 1$ cannot occur under the hypotheses of Theorem~\ref{local_thm}, neither can $F_{2,\bar{w}}\oplus \ZZ_p\oplus A$ nor $\cyclic{p}\oplus \ZZ_p\oplus A$.

\begin{proof}
 Let $M=E(K)\hatox\ZZ_p$.
 Since the formal logarithm induces a $\QQ_p[G]$-isomorphism $M\otimes\QQ_p\to K$, we have $M\otimes\QQ_p\cong \QQ_p[G]$.
 Since $K$ is unramified, the $p$-th roots of $1$ cannot be contained in $K$ and therefore $M_t= E(K)[p^{\infty}]$ is cyclic. Hence the classification in Theorem~\ref{cyclic_class_thm} applies.

 From Lemma~\ref{unram_lem}, we know that $D = \hat H^0(G,M)$ is cyclic of order $p$ if the reduction is split multiplicative and $p\mid c_v$ and otherwise it is trivial.

 Assume $D$ is trivial.
 If the reduction is split multiplicative then $p\nmid c_v$ and as the Tate parameter $q$ of $E$ has a valuation that is not a multiple of $p$, there cannot be any points of order $p$ on $E(\QQ_p)$.
 For other types of reduction, we know that both the formal group $\hat E(p\ZZ_p)$ and the group of components $\Phi(\FF_p)$ have no elements of order $p$, which means that $E(\QQ_p)[p^\infty]$ injects into $\tilde{\mathcal{E}}_v^0(\FF_p)$.
 By the Hasse-Weil bound in the case of good reduction and by direct considerations in the case of bad reduction, we deduce that $E(\QQ_p)[p^{\infty}]$ is either trivial or cyclic of order $p$.
 In the first case $M_t=E(K)[p^\infty]$ is also trivial and the only option for $M$ is then $\ZZ_p[G]$.
 In the second case, $M_t$ is either $\cyclic{p}$ or $F_{2,\bar{w}}$, which explains the other two entries in the first table.

 We can now assume that the reduction is split multiplicative and $p$ divides $c_v$.
 So $D$ is cyclic of order $p$.
 Set $j=\ord_p(c_v)>0$ and $i$ such that $p^i$ is the order of $E(\QQ_p)[p^\infty]$.

 Since $M_t$ injects into $\Phi(\FF_w)[p^\infty]\cong \cyclic{p^j}$, we must have $M_t\cong\cyclic{p^i}$ as a $\ZZ_p[G]$-module.
  Therefore, we must have $j\geq i$.
 If $i=0$, then the classification limits us to only one option, namely $M\cong\ZZ_p\oplus A$.
 If $i>0$, then let $P$ be a point of exact order $p^i$ on $E(\QQ_p)$ and consider the isogeny from $E$ to $E'$ whose kernel is generated by $P$.
 In terms of the Tate curve, this map $K^\times/q^\ZZ \to K^\times/q'^{\ZZ}$ is induced from the identity map and $(q')^{p^i} = q$.
 Now the order of the torsion subgroup of $E'$ is no longer divisible by $p$ as $q'$ is not a $p$-th power in $\QQ_p$.
 This implies that $E'(K)\hatox \ZZ_p$ is isomorphic to $\ZZ_p\oplus A$ or $\ZZ_p[G]$ depending on whether the Tamagawa number $c'_v$ of $E'$ is divisible by $p$ or not.
 If $p\nmid c'_v$, then $j=i$ and the extension between $\cyclic{p^i}$ and $\ZZ_p[G]$ must split.

 Instead if $p\mid c'_v$, then $j>i$.
 The only extensions of $\cyclic{p^i}$ and $\ZZ_p\oplus A$ with $D\cong\cyclic{p}$ are $\bigl\{\cyclic{p^i}\bigm\vert \ZZ_p\bigr\}\oplus A$ or $\bigl\{\cyclic{p^i}\bigm\vert A\bigr\}\oplus \ZZ_p$.
 We can exclude the first case because the map from $M^G=E(\QQ_p)\hatox \ZZ_p$ to $E'(\QQ_p)\hatox\ZZ_p$
 must be surjective.
\end{proof}

The same proof should work if $k=\QQ_p$ is replaced by a finite extension of $\QQ_p$ with ramification index less than $p-1$.
The case of reduction IV and IV${}^*$ when $p=3$ can be treated as well but they are more complicated as illustrated by the last example in this section.

\subsection{Examples}
We wish to give explicit examples for Theorem~\ref{local_thm} and thus we take $K/\QQ_p$ to be the unramified extension of degree~$p$.
To cover all possible cases with split multiplicative reduction it is enough to look at Tate curves whose parameters are, say,
\[ q = p^{p^j}\cdot (1+p)^{p^i} = \Bigl(p^{p^{j-i}}\cdot (1+p)\Bigr)^{p^i} \]
with integers $j\geq i \geq0$.
This $q$ is a $p^i$-th power so there are $p^i$-torsion points in $E(\QQ_p)$.
Since $1+p$ is not a $p$-th power, we get $E(\QQ_p)[p^{\infty}]=\cyclic{p^i}$.
Also $v_p(c)= v_p(v_p(q)) =j$.

For the two additive cases, we can take the following two examples.
First $E:y^2+y=x^3$ over $\QQ_3$ has additive reduction of type II.
There is a $3$-torsion point $(0,0)$ in $E(\QQ_3)$.
Hence $E(K)\hatox\ZZ_p$ must be $\bigl\{\cyclic{p}\bigm\vert\ZZ_p\oplus A\bigr\}$.

Secondly, $E:y^2+y=x^3-270\,x-1708$ has additive reduction of type II${}^*$ over $\QQ_3$.
This time one can show that $E(\QQ_3)$ does not contain a point of order $3$ directly by checking the roots of the $3$-division polynomial.
Alternatively on can use the map $\partial$ in the exact sequence
\[ \xymatrix@1{ 0\ar[r] & E(\QQ_p)[p]\ar[r] & \tilde E^0(\FF_p)[p] \ar[r]^(0.4){\partial} & \hat E(p\ZZ_p)/p \hat E(p\ZZ_p).}\]
For instance the point $P=(1,3+3^2+2\cdot 3^3+\cdots)$ has non-trivial non-singular reduction.
Then $Q=3P = (3^{-2}+2\cdot 3^2+\cdots, 2\cdot 3^{-3} + 2\cdot 3^{-2}+\cdots)$ belongs to $\hat E(p\ZZ_p)$ but not to $p \hat E(p\ZZ_p)= \hat E(p^2\ZZ_p)$.
Therefore $\partial(\tilde P) =Q + p \hat E(p\ZZ_p)$ is non-trivial.
It follows that $E(\QQ_p)[p]$ is trivial and $E(K)\hatox \ZZ_p$ is isomorphic to $\ZZ_p[G]$.

Now to curves with good reduction.
The curve $y^2 + x\,y + y = x^3 - 171\,x - 874$ has good ordinary, anomalous reduction $\tilde P = (1,0)$ over $\QQ_3$.
With the methods from the previous case one can show that $E(\QQ_p)[p]$ is trivial in this case.
Thus $E(K)\hatox\ZZ_p$ is a free $\ZZ_p[G]$-module.

Instead the curve $E:y^2 + x\,y + y = x^3 + 4\,x - 6$ has also good ordinary, anomalous reduction over $\QQ_3$, but it contains a $3$-torsion point $(2,2)$.
It can be shown that this point does not become divisible by $3$ in $K$.
Therefore $E(K)\hatox\ZZ_p$ must be isomorphic to $\bigl\{\cyclic{p}\bigm\vert \ZZ_p\oplus A\bigr\}$.

Finally, here an example when $E(K)\hatox\ZZ_p$ must be $F_{2,\bar{w}}\oplus \ZZ_p [G]$.
The curve $E:y^2+y=x^3+x^2+x$ has a $3$-torsion point $T=(0,0)$ and it has good reduction over $\QQ_3$.
The extension $K$ can be obtained by adjoining $t$ to $\QQ_3$ with $t^3+2\,t+1=0$.
Then the point
\begin{multline*}
  S = \bigl((t^2 + 2\,t + 2) + (2\,t^2 + t + 1)\cdot 3 + (t^2 + t + 2)\cdot 3^2 + (2\,t + 2)\cdot 3^3 + \cdots,\\
  (t^2 + t + 1) + (t + 2)\cdot 3 + t^2\cdot 3^2 +  (2\,t^2 + 2\,t + 1)\cdot 3^3 +\cdots \bigr)
\end{multline*}
is such that $3\,S = T$.

Therefore we have found explicit examples for all possible $\ZZ_p G$-modules in Theorem~\ref{local_thm}.

To illustrate that the situation is more complicated in the cases not treated in Theorem~\ref{local_thm}, we add one example.
The curve
\[ y^2 + xy + 9\, y = x^3 -x^2+9\,x+9 \]
over $\QQ_3$ has additive reduction of type IV with Tamagawa number $c_v=3$.
This curve has a rational $9$-torsion point $T$ with $x$-coordinate $3 + 2\cdot 3^2 + 2\cdot 3^4 + 3^5 + 3^6 + 2\cdot 3^7 + 3^8 + 2\cdot 3^9 + \mathbf{O}(3^{10})$. 
It has bad reduction, but $3T$ is a $3$-torsion point with good reduction.
It seems that $E(K)\hatox \ZZ_p$ is isomorphic to $\bigl\{\cyclic{9}\bigm\vert \ZZ_3\bigr\} \oplus A$.

\section{Descent for Mordell-Weil and Selmer groups}\label{descent_sec}

We pass to studying global extensions and gather the tools to investigate how the Galois group of a finite extension acts on the Mordell-Weil group and on the Selmer group.

Let $k$ now be a number field and let $K/k$ be a finite extension with Galois group $G$.
Let $E/k$ be an elliptic curve and let $p$ be an odd prime.

We write $M=E(K)\hatox \ZZ_p=E(K)\otimes \ZZ_p$, which we are going to study as a $\ZZ_p[G]$-module.
Recall that $M_t$ is the torsion subgroup of $M$ and $M_f=M/M_t$.
We have $M^G = E(k)\otimes \ZZ_p$.
Similarly the $G$-fixed part of $M\otimes \QQ_p = E(K)\otimes \QQ_p$ is $E(k)\otimes \QQ_p$.
Consider instead the limit $\varinjlim E(K)/p^n E(K)$ which naturally identifies with $E(K)\otimes \QZp$.
The map comparing the $G$-fixed part of $E(K)\otimes \QZp$ with $E(k)\otimes \QZp$ measures if any point $P\in E(k)$ becomes divisible by $p$ in $E(K)$ when they were not in $E(k)$:
\begin{lem}\label{div_lem}
  We have an exact sequence
  \begin{equation*}
   \mathclap{ \xymatrix@1@C-1ex{
      0 \ar[r] & \ker \Bigl(H^1(G,M_t)\to H^1(G,M)\Bigr) \ar[r] & E(k)\otimes \QZp \ar[r] & \Bigl(E(K)\otimes \QZp\Bigr)^G \ar[r] & H^1(G, M_f)\ar[r] & 0.
    }}
  \end{equation*}
\end{lem}
\begin{proof}
  First the definition of $M_f$ yields the long exact sequence
  \[
    \xymatrix@1{
      0 \ar[r] & E(k)[p^{\infty}] \ar[r] & E(k)\otimes \ZZ_p \ar[r] & M_f^G\ar[r] & H^1(G,M_t) \ar[r] & H^1(G,M).
    }
  \]
  Further we have an exact sequence
  \[
    \xymatrix@1{
    0\ar[r] & M_f^G \ar[r] & E(k)\otimes \QQ_p \ar[r] & \Bigl(E(K)\otimes \QZp\Bigr)^G \ar[r] & H^1(G,M_f)\ar[r] & 0
    }
  \]
  deduced from the short exact sequence
  \[
    \xymatrix@1{0\ar[r] & M_f \ar[r] & E(K)\otimes \QQ_p \ar[r] & E(K)\otimes \QZp \ar[r] &0 }
  \]
  and an isomorphism
  \begin{equation}\label{h2mf_eq}
    H^1\bigl(G, E(K)\otimes \QZp\bigr) \cong H^2(G,M_f)
  \end{equation}
  which will be useful later.
  The kernel-cokernel sequence for the composition $E(k)\hatox \ZZ_p \to E(k)\otimes \QQ_p$ via $M_f^G$ produces the exact sequence in the lemma.
\end{proof}

In particular, if $M$ is torsion-free, then $H^1(G,M)$ is the cokernel of the injective map $E(k)\otimes \QZp\to \bigl(E(K)\otimes\QZp\bigr)^G$.

We can compare $M$ to the local group of points.
Let $v$ be a finite place of $k$ and let $w$ be a place of $K$ above $v$.
Let $M_w = E(K_w)\hatox \ZZ_p$.
If $v\nmid p$ then $M_w$ is finite; otherwise it is a $\ZZ_p$-module of rank $[K_w:\QQ_p]$.
In Section~\ref{local_mw_sec} we have investigated $M_w$ as a $\ZZ_p\bigl[\Gal(K_w/k_v)\bigr]$-module.
Setting $M_v = \bigoplus_{w\mid v} M_w$ we get a $\ZZ_p[G]$-module.

By $S$ we will denote a fixed choice of a finite set of places in $k$ containing all places above $p$ and $\infty$, such that $E$ has good reduction outside $S$ and such that $K/k$ is unramified outside~$S$.
The Galois group of the maximal extension of $k$ which is unramified outside $S$ is denoted by $G_S(k)$.

We define the $p$-primary Selmer group $\Sel_k$ as the subgroup of $H^1\bigl(G_S(k), E[p^{\infty}]\bigr)$ consisting of elements whose restrictions at all places $v\in S$ lie in the image of the local Kummer map $\kappa_v \colon E(k_v)\otimes \QZp \to H^1\bigl(k_v, E[p^{\infty}]\bigr)$.
The compact $p$-adic Selmer group $\Selc_k$ is the subgroup of $H^1\bigl(G_S(k),T_pE\bigr)$ defined in the same way with respect to the local Kummer maps $E(k_v)\hatox \ZZ_p \to H^1\bigl(k_v, T_pE\bigr)$ in the cohomology of the Tate module~$T_pE$.
The cokernel of the global Kummer map $E(k)\otimes \QZp\to \Sel_k$ is the $p$-primary part of the Tate-Shafarevich group $\Sha(E/k)$.
We will simply denote it by $\Sha_k$.

\begin{assu}\label{sha_assu}
  Throughout this article we assume that the $p$-primary subgroup $\Sha_k$ of the Tate-Shafarevich group is a finite group for all number fields $k$ and all elliptic curves.
\end{assu}
One can drop this assumption and work with the quotient of $\Sha_k$ by its maximal divisible subgroup instead.
It follows from this assumption that $\Selc_K$ is simply $M = E(K)\otimes \ZZ_p$ and $\Selc_k$ is~$M^G$.

Using the isomorphism~\eqref{h2mf_eq}, we obtain a commutative diagram
\[
  \xymatrix@R-2ex{
   0 \ar[r] & \Bigl(E(K)\otimes \QZp\Bigr)^G \ar[r] & \Sel_K^G \ar[r] & \Sha_K^G \ar[r] & H^2\bigl(G,M_f\bigr)\\
   0 \ar[r] & E(k)\otimes \QZp \ar[r] \ar[u] & \Sel_k \ar[r] \ar[u]_{\alpha} & \Sha_k\ar[r]\ar[u]_{\eta} & 0
  }
\]
which describes the link of the exact sequence in Lemma~\ref{div_lem} to the natural restriction maps $\alpha$ on the Selmer group and $\eta$ on the Tate-Shafarevich groups.
Below we will determine the kernel and cokernel of $\alpha$.
We denote the kernel of $\eta$ by $C_{K/k}$ and call it the capitulation subgroup of $\Sha(E/k)[p^\infty]$ with respect to the extension $K/k$.
The snake lemma applied to this diagram provides us with an exact sequence
\begin{equation}\label{shadescent_seq}
  \xymatrix{
    & \ker(\alpha) \ar[r] & C_{K/k} \ar[r] & H^1(G, M_f) \ar `r[d] `[ll] `l[dlll] `d[dll] [dll] & \\
    & \coker(\alpha) \ar[r] & \coker(\eta) \ar[r] & H^2(G, M_f).
  }
\end{equation}
\begin{assu}\label{tors_assu}
  We assume that $E(K)$ does not contain any non-trivial $p$-torsion elements.
\end{assu}

\begin{lem}\label{sha_descent_lem}
  Under Assumption~\ref{tors_assu}, we have the exact sequence
  \[\xymatrix@1{
    0 \ar[r] &\ker(\alpha)\ar[r] & C_{K/k} \ar[r] & H^1(G, M) \ar[r] &  \coker(\alpha) \ar[r] & \coker(\eta) \ar[r] & H^2(G, M).
  }\]
\end{lem}
\begin{proof}
  The kernel of the first map in the sequence~\eqref{shadescent_seq} is in fact the kernel of $H^1(G,M_t)\to H^1(G,M)$ as we have seen in Lemma~\ref{div_lem}.
  By the assumption $M$ is $\ZZ_p$-free, hence $M_t$ is trivial and $M_f=M$.
\end{proof}

Our next aim is to determine the cokernel of the restriction map $\alpha$ between $\Sel_k$ and $\Sel_K^G$.
Global duality, in form of the exact sequence due to Cassels, see Section~1.7 in~\cite{rubineulersystems}, allows us to compare $M$ with the local terms $M_v$.
It yields the exact sequence
\[
  \xymatrix@1{M\ar[r] & \bigoplus_{v\in S} M_v \ar[r] & H^1\bigl(G_S(K),E[p^\infty]\bigr)^\vee \ar[r] & \Sel_K^{\vee} \ar[r] & 0}
\]
of $\ZZ_p[G]$-modules.
The kernel of the first map is defined to be the fine (or strict) Mordell-Weil group $\mathfrak \Relc_{S,K}$.
In many circumstances this group is trivial; see for instance~\cite{waldschmidt}.
By global duality, as in Theorem~8.6.8 in~\cite{cnf}, $\Relc_{S,K}$ is dual to $H^2\bigl(G_S(K), E[p^\infty]\bigr)$.
We compare the $G$-fixed part of the dual of the above exact sequence with the corresponding sequence over $k$ to make the map $\alpha$ appear in the following large commutative diagram.
The top sequence is only a complex; at the terms where the complex is not necessarily exact we use the symbol $\xymatrix@1{{}\ar@{o->}[r]&{}}$.
\[
\mathclap{
  \xymatrix@R-1ex{
    0\ar[r] & \Sel_K^G \ar[r] & H^1\bigl(G_S(K),E[p^{\infty}]\bigr)^G \ar[r] & \bigoplus_{v\mid S} \bigl(M_v^{\vee}\bigr)^G \ar@{o->}[r] & \bigl(M^{\vee}\bigr)^G\ar@{o->}[r] & \bigl(\Relc_{S,K}^{\vee}\bigr)^G \ar@{o->}[r] & 0 \\
    0\ar[r] & \Sel_k \ar[r]\ar[u]^\alpha & H^1\bigl(G_S(k),E[p^{\infty}]\bigl) \ar[r]\ar[u]^{\beta} & \bigoplus_{v \in S} \bigl( M_v^G\bigr)^\vee \ar[r]\ar[u]^{\gamma} & \bigl(M^G\bigr)^{\vee}\ar[r]\ar[u]_{\delta} & \Relc_{S,k}^{\vee} \ar[r]\ar[u]_{\varepsilon} & 0
  }}
\]
The map $\gamma$ is then $\oplus_{v\in S} \gamma_v$ where $\gamma_v$ is the dual of the natural norm map; the kernel of $\gamma_v$ is dual to
\begin{align}\label{Dv_eq}
  \begin{split}
    D_v &= \coker \Bigl(\No\colon \bigoplus_{w\mid v }E(K_w)\hatox \ZZ_p\to E(k_v)\hatox \ZZ_p\Bigr) \\
    & = \hat H^0\Bigl( G, \bigoplus_{w \mid v} E(K_w)\hatox\ZZ_p\Bigr) \cong \hat H^0\bigl(G_w, M_w\bigr).
  \end{split}
\end{align}
for any chosen $w$ above $v$.
Here $G_w = \Gal(K_w/k_v)$ is the decomposition group.
The groups $D_v$ were studied in detail in the Section~\ref{local_norm_sec}.

The maps $\delta$ and $\varepsilon$ are the dual of the norm map
\[
 \hat\delta\colon\bigl(E(K)\otimes \ZZ_p\bigr)_G \to E(k)\otimes\ZZ_p
\]
and its restriction to the fine Mordell-Weil group.
Therefore $\ker(\delta)$ is dual to $\hat H^0(G,M)$ and $\coker(\delta)$ is dual to $\hat H^{-1}(G,M)$.
The commutativity of the diagram follows from the functoriality of global duality and the local duality of restriction and corestriction.

Our assumption~\ref{tors_assu}, that $E(K)$ has no non-trivial $p$-torsion points, implies that $\beta$ is an isomorphism by the inflation-restriction-transgression exact sequence, see Proposition~1.6.6 in~\cite{cnf}.

We are in a situation where we have a morphism of complexes $A^\bullet \to B^\bullet$ with $A^\bullet$ exact.
Let $\bar A^\bullet$ be the complex of kernels and $\bar B^\bullet$ the complex of cokernels.
As a consequence of the long exact sequences of cohomology of complexes in short exact sequences, we can deduce that there is a long exact sequence
\[
 \xymatrix@1{
  \cdots \ar[r] & H^{i+1}(\bar A^\bullet) \ar[r] & H^i(B^\bullet)\ar[r] & H^i(\bar B^\bullet)\ar[r] & H^{i+2}(\bar A^\bullet) \ar[r] & \cdots
  }
\]
In our case, since the first two terms of $B^\bullet$ have trivial cohomology, we deduce that $\ker(\alpha)=0$ and that
\[
  \xymatrix@1{
    0\ar[r] & \coker(\alpha)\ar[r] & \ker(\gamma) \ar[r] & \ker\Bigl( \ker\delta\to \ker\varepsilon \Bigr)\ar[r] & \cdots
  }
\]
is exact.
The image at the end is a subquotient of $\bigoplus_{w\in S_K} \bigl(E(K_v)\hatox\ZZ_p\bigr)^\vee$.

If we are just interested in $\coker(\alpha)$ we see it here as the kernel of $\ker\gamma\to \ker\delta$.
The dual of this map is $\hat H^0(G,M) \to \bigoplus_{v\in S} D_v$.
As $E(k)\otimes\ZZ_p\to\hat H^0(G,M)$ is surjective, we find that $\coker(\alpha)$ is dual to the cokernel of $E(k)\otimes\ZZ_p\to \bigoplus_{v \in S} D_v$.
We have shown the following:

\begin{prop}\label{cokeralpha_prop}
  We have $\ker(\alpha)=0$ and $\coker(\alpha)$ is dual to the cokernel of $E(k)\otimes \ZZ_p \to D_{K/k}$ where $D_{K/k} = \bigoplus_{v\in S} D_v = \bigoplus_{v \in S} \hat H^0\bigl(G_w,E(K_w)\hatox\ZZ_p\bigr)$.
\end{prop}

We have calculated $D_v$ in many circumstance in Section~\ref{local_norm_sec}.
Here is a summary for the most frequent situation

\begin{prop}\label{easy_d_prop}
  Suppose that no place above $p$ is wildly ramified in $K/k$ and suppose that the ramification index $e_v$ is not divisible by $p$ at all places where $E$ has bad reduction.
  Then
  \[
    D_{K/k} \cong \bigoplus_{v\in S_{\mathrm{b}}} \ZZ/(c_v,f_v,p^\infty)\ZZ \oplus \bigoplus_{v\in S_{\mathrm{r}}} \tilde E(\FF_v)[p^\infty]/e_v \tilde E(\FF_v)[p^\infty]
  \]
  where $S_\mathrm{b}$ is the set of all places where $E$ has bad reduction and $S_{\mathrm{r}}$ is the set of all places where $p$ divides the ramification index $e_v$.
\end{prop}

\begin{proof}
  For all places with $p\nmid e_v$, Proposition~\ref{notwild_prop} shows that $D_v$ is cyclic of order equal to $\gcd(c_v,f_v,p^\infty)$.
  This is non-trivial only for places of bad reduction and they appear in the first sum above.

  If $v$ is a place in $S_{\mathrm{r}}$, then, by assumption, $E$ has good reduction and $v\nmid p$.
  Therefore Proposition~\ref{good_not_p_prop} applies and gives the second term.
\end{proof}

Note that the map $E(k)\to E(k_v) \to D_v$ is explicit and easy to calculate.
As a consequence, we can effectively determine the cokernel of $\alpha$ using Proposition~\ref{cokeralpha_prop}.
This is efficient, but more involved if the conditions of Proposition~\ref{easy_d_prop} do not hold, but local ad hoc calculations will allow to determine $D_{K/k}$ and hence $\coker(\alpha)$.
As a consequence, the above proposition can be used to determine $\coker(\alpha)$ explicitly using only local information and information about $E$ over $k$ in all cases.

\section{The Galois module structure of Mordell-Weil groups}\label{galmod_sec}

Let $E/k$ be an elliptic curve and let $K/k$ be a Galois extension of number fields with group $G$.
Let $p$ be an odd prime.
Throughout this section, we continue to work under the Assumptions~\ref{sha_assu} and~\ref{tors_assu}.

Recall that the aim is to understand in what cases we can determine the structure of $M=E(K)\hatox\ZZ_p$ as a $\ZZ_p[G]$-module, preferably with information that is easier to access than computing $M$ itself.

We first recall the results of~\cite{bmw1}, which can be extended to all situations when we have ``perfect control'', i.e., when
\[\alpha\colon \Sel_k \to \Sel_K^G \]
is an isomorphism.

\begin{thm}[Theorem~2.2 in~\cite{bmw1}]
  Assume that $\coker(\alpha)$ is trivial.
  Fix a $p$-Sylow subgroup $\mathcal{H}$ of $G$.
  Then $M$ is a projective $\ZZ_p[G]$-module if and only if $M\otimes\QQ_p$ is a free $\QQ_p[\mathcal{H}]$-module and $C_{K/K^H}$ is trivial for all subgroups $H$ in $\mathcal{H}$.
\end{thm}

When the $p$-Sylow subgroup $\mathcal{H}$ is cyclic, we can use the results of Yakovlev~\cite{yakovlev} to determine $M$ as follows.
Since $\mathcal{H}$ is cyclic, say of order $p^n$, we have a tower of fields $F_i$ such that $[K:F_i]=p^{n-i}$ and $F_0 = K^{\mathcal{H}}$.

\begin{thm}[Theorem~2.6 in~\cite{bmw1}]
  Suppose that the $p$-Sylow subgroup $\mathcal{H}$ of $G$ is cyclic.
  Assume that $\coker(\alpha)$ is trivial.
  Then $M$ is determined up to isomorphism by the the ranks of $E(F_i)$ and the knowledge of the capitulation kernels as $\ZZ_p\bigl[N_G(\mathcal{H})\bigr]$-modules together with the restriction and corestriction maps between them:
  \[
    \xymatrix@1{C_{K/F_0} \ar@<.5ex>[r] & \ar@<.5ex>[l] C_{K/F_1} \ar@<.5ex>[r] & \ar@<.5ex>[l]\cdots \ar@<.5ex>[r] & \ar@<.5ex>[l] C_{K/F_{s-1}}}.
  \]
\end{thm}

Since $C_{K/F_i} \cong \hat H^{-1}\bigl(K/F_i, M\bigr)\approx H^1\bigl(K/F_i, M\bigr)$ the diagram of restrictions and corestrictions above is an example of a Yakovlev diagram as defined in Section~6.2 in~\cite{torzewski}, where the precise formulation of this theorem is discussed.
Specialising to the situation when the $p$-Sylow $\mathcal{H}$ of $G$ is cyclic of order $p$, the Yakovlev diagram simplifies then to a single group $H^1\bigl(\mathcal{H},M)$ viewed as a $\FF_p\bigl[N_G(\mathcal{H})/\mathcal{H}\bigr]$-module.
It is well possible that the results in~\cite{macias17} also hold under the weaker hypothesis that $\alpha$ is surjective.

In order to generalise to the case when $\alpha$ is not necessarily surjective, we will use Torzewski's result~\cite{torzewski} which extends Yakovlev's theorem (as in Theorem~\ref{torzewski_thm} for the dihedral case) by involving the regulator constants $s_{\Theta}(M)$, too.
This new ingredient can be linked to arithmetic information as follows.
Fix an invariant differential $\omega$ on $E$ and write $u_v = \lvert \omega/\omega_{\text{Néron}} \rvert_v$ for when it differs from the Néron differential $\omega_{\text{Néron}}$ of $E$ at the finite place $v$.
For any field $F/k$, we define $C(E/F) = \prod_v c_v(E/F) \cdot u_v$ to be the modified product over all finite places $v$ of $F$ of the Tamagawa numbers $c_v$.
This quantity together with the real and complex periods with respect to $\omega$ should appear in the leading term of the Birch and Swinnerton-Dyer conjecture for $E/F$.
The following is a consequence of Theorem~2.3 in~\cite{DD_square}.

\begin{thm}\label{bartel_thm}
   Assume that no place of additive reduction ramifies in $K/k$.
   Let $\Theta = \sum m_i\, H_i$ be a Brauer relation for $G$.
   For $F= K^{H_i}$ write $s_i = \ord_p(\lvert\Sha_{F}\rvert) + \ord_p\bigl( C(E/F) \bigr)$.
   Then $s_{\Theta} (M) = - \sum_{i} m_i s_i$.
\end{thm}

Under our assumption, we can determine the parity of $s_{\Theta}(M)$ from the Tamagawa numbers only as the Tate-Shafarevich groups are of square order.

\begin{cor}\label{parity_cor}
  Under our assumptions, $s_{\Theta}(M) \equiv \sum_i m_i\, \ord_p\bigl(C(E/F^{H_i})\bigr)\pmod{2}$ for any Brauer relation $\Theta = \sum_i m_i \, H_i$.
\end{cor}

The parity here will link directly to global root numbers.
However in the general case, we often need the ranks rather than just their parity.
In the particular cases that we turn our attention to, the valuation of $s_{\Theta}(M)$ and the rank over $k$ will provide more information than the root numbers.

\subsection{Cyclic extensions}
Suppose first that $G$ is a cyclic group of order $p$.
Recall that we write $D_{K/k}$ for $\bigoplus_{v\in S} D_v$, which is a $\FF_p$-vector space in our situation.

\begin{cor}\label{cyclic_d_cor}
  Let $S_{\mathrm{r}}^0$ be the set of all ramified places in $k$ not lying above $p$ and at which $E$ has good reduction.
  Then
  \[ \dim_{\FF_p} D_{K/k} \geq \#\bigl\{ v \bigm\vert v\text{ inert  and } p\mid c_v\bigr\} + \sum_{v \in S_{\mathrm{r}}^0} \dim_{\FF_p} \tilde E(\FF_v)[p].\]
  If no place of bad reduction and no place above $p$ ramifies, then we have equality.
\end{cor}

\begin{proof}
  First the dimension of $D_{K/k}$ is larger than the sum of the dimensions of $D_v$ for all places excluding the ramified places above $p$ and the ramified places at which $E$ has bad reduction.
  For the remaining $v$, we calculated $D_v$ in Proposition~\ref{easy_d_prop} and we can simplify it a bit because $G$ is cyclic of order $p$.
\end{proof}

Recall that, by Proposition~\ref{cyclic_lattices_prop}, there are only three indecomposable $\ZZ_p[G]$-lattices $\ZZ_p$, $A$ and $\ZZ_p[G]$ in this case.
We note that $M$ cannot be determined by $\rk E(K)$ and $\rk E(k)$ only, but they will determine it together with the order of $H^1(G,M)$.

\begin{prop}\label{cyclic_items_prop}
  Let $E/k$ be an elliptic curve and $K/k$ a cyclic extension of degree $p$. 
  \begin{enumerate}[nosep,label={(\roman*)}]
    \item If $\Sha_k$ and $D_{K/k}$ are trivial, then $\rk E(K)$ and $\rk E(k)$ determine $M$.
    \item \label{pp2_item} If $\rk E(k)=0$ and $\Sha_k$ is trivial, then $\rk E(K) \leq (p-1) \dim_{\FF_p} D_{K/k}$.
    \item \label{pp3_item} If $\rk E(k)=0$ and $\Sha_k$ is trivial, but  $\rk E(K)< (p-1) \dim_{\FF_p} D_{K/k}$, then $\Sha_K$ is not trivial.
    \item \label{pp4_item} If $\alpha$ is surjective, $\rk E(K)> \rk E(k) = 1$, and $\Sha_k$ trivial, then $M\cong \ZZ_p[G]$ and $\Sha_K=0$.
  \end{enumerate}
\end{prop}

\begin{proof}
  For the first point, the hypothesis imply that $\coker(\alpha)$ and $C_{K/k}$ are trivial.
  By Lemma~\ref{sha_descent_lem}, this implies that $H^1(G,M)=0$ and hence $M$ is a direct sum of copies of $\ZZ_p$ and $\ZZ_p[G]$, which can be determined by the ranks alone.

  In part~\ref{pp2_item} and~\ref{pp3_item}, $\rk E(k)=0$ implies that $M$ is a a direct sum of copies of~$A$.
  The number of copies is equal to $\dim_{\FF_p} H^1(G,M)$ which is the dimension of $\ker \bigl(\coker(\alpha) \to \coker(\eta)\bigr)$.
  As the rank is zero over $k$, the cokernel of $\alpha$ is dual to $D_{K/k}$.
  If the resulting inequality is strict, then $\coker(\eta)$ is non-trivial and hence so is $\Sha_K$.

  For the final part~\ref{pp4_item}, the surjectivity of $\alpha$ and the triviality of $\Sha_k$ imply that $H^1(G,M)=0$.
  As the rank grows but is equal to $1$ over $k$, we must have $M=\ZZ_p[G]$.
  This now also implies that $H^2(G,M)=0$.
  From Lemma~\ref{sha_descent_lem} we learn that $\eta$ is surjective.
  Since $G$ is a $p$-group, the triviality of $\Sha_K^G$ implies that $\Sha_K$ is trivial.
\end{proof}

\begin{proof}[Proof of Theorem~\ref{ex2_thm}]
  The assumption that $L(E,\chi,1)\neq 0$ implies by Kato's result in~\cite{kato} that $\rk E(K) = \rk E(k)$.
  Therefore $M\cong\ZZ_p^r$ with $r=r_\QQ$.
  By Corollary~\ref{cyclic_d_cor}, the dimension of $D_{K/\QQ}$ is greater or equal to $u_1+u_2$ because the sum of the dimensions in the second term is larger or equal to the number of $\ell \in S_{\mathrm{r}}^0$ such that $p\mid \#\tilde E(\FF_{\ell})$, which equals $u_2$.
  Therefore,  Proposition~\ref{cokeralpha_prop} tells us that $\coker(\alpha)$ has dimension at least equal to $u_1+u_2-r$.
  Since $H^1(G,M)=0$, $\coker(\alpha)$ injects into $\coker(\eta)$ by Lemma~\ref{sha_descent_lem}.
  Hence $\dim_{\FF_p} \Sha_K[p] \geq \dim_{\FF_p}\Sha_K^G[p]$, which is the minimal number of generators of $\Sha_K^G$, and it must be equal or larger than $\dim_{\FF_p}\coker(\eta) \geq u_1+u_2-r$.
  (By the way, we also get $\dim_{\FF_p} \Sha_K^G[p] \leq u_1+u_2 +r + \dim_{\FF_p} \Sha_k[p]$.)
\end{proof}

The exact sequence relating the Mordell-Weil group $M$ and the Tate-Shafarevich group $\Sha_K$ to the Selmer group $\Sel_K$ is split as a sequence of abelian groups~\cite{gillibert}, but not necessarily as $G$-modules.
There is one important exception.
\begin{lem}
  Suppose $M\cong \ZZ_p[G]$.
  Then \[\Sel_K \cong \bigl(E(K)\otimes\QZp\bigr) \oplus \Sha_K\] as $G$-modules and we have an exact sequence
  \[
    \xymatrix@1{0\ar[r] & \Sha_{k}\ar[r] & \Sha_K^G \ar[r] & D_{K/k}^{\vee}\ar[r] & 0. }
  \]
\end{lem}
\begin{proof}
  Note that there is a $G$-equivariant morphism $\Hom_{\ZZ_p}\bigl(M,\ZZ_p\bigr) \to \Hom_{\ZZ_p}\bigl( M\otimes \QZp,\QZp\bigr)$ sending $f$ to the map $m\otimes x \mapsto f(m)\cdot x$ for all $m\in M$ and~$x\in\QZp$.
  By the assumption that $E(K)[p]$ is trivial, this is an isomorphism.
  If $M$ is $\ZZ_p[G]$-free, then so is the Pontryagin dual of~$E(K)\otimes\QZp= M\otimes \QZp$.
  Since the quotient of $\Sel_K^{\vee}$ by $\Sha_K^{\vee}$ is now a projective $\ZZ_p[G]$-module, we must have $\Sel_K^{\vee}\cong\Hom\bigl(M,\ZZ_p\bigr)\oplus\Sha_K^{\vee}$ whose Pontryagin dual is the initial statement.
  The exact sequence~\eqref{shadescent_seq}, together with $H^1(G,M)=0$ and $H^2(G,M)=0$, shows $C_{K/k}=0$ and that $\coker(\eta)\cong\coker(\alpha)$.
  As the norm map $M_G\to M^G$ is an isomorphism, the cokernel of $\alpha$ is dual to~$D_{K/k}$.
\end{proof}

\begin{thm}\label{cyclic_pos_thm}
  Let $E/k$ be an elliptic curve and $p$ an odd prime.
  Suppose that $\rk E(k)=0$, that $\Sha_k=0$ and that the image of the Galois representation $\Gal(\bar{k}/k)\to \GL\bigl(E[p]\bigr)$ contains $\SL\bigl(E[p]\bigr)$.
  Then there is a positive proportion of cyclic extensions $K/k$ of degree $p$ with prime conductor such that $M=0$.
\end{thm}

Note that we expect that a positive proportion of elliptic curves $E/k$ should satisfy the hypothesis in the theorem.
It is important to emphasis that this results is weaker than Theorem~9.21 in~\cite{mazur_rubin_ds} in the sense that we restrict to curves of rank $0$ and prime degree, but more general in other aspects.
Our proof imitates closely their work and one could also use their Section~11 to obtain a quantitative result.

\begin{proof}
  Fix $\Sigma$ to be a finite set of places containing all bad places $v$ in $k$, all places such that $p\mid c_v$, and assume it is sufficiently large so that $H^1\bigl(G_{\Sigma}(k), \mu[p]\bigr) = 0$.
  Among all prime ideals outside of $\Sigma$ and outside the bad places for $E$, we let $\mathcal{P}$ be the subset of prime ideals $\mathfrak p$ verifying the following two conditions: Firstly, that $N(\mathfrak{p}) \equiv 1 \pmod{p}$ and, secondly, that all $\Sigma$-units are $p$-th powers in $\mathcal{O}_{\mathfrak{p}}^{\times}$.
  For each $\mathfrak p$ in~$\mathcal{P}$, there is a Galois extension $k(\mathfrak{p})/k$ of degree $p$ exactly ramified at $\mathfrak{p}$ such that all places in $\Sigma$ are split; this follows from the global class field argument in the proof of Lemma~9.15 in~\cite{mazur_rubin_ds}.
  Let $\mathcal{P}_0$ be the subset of primes in $\mathcal{P}$ such that $p\nmid \# \tilde E(\FF_{\mathfrak{p}})$.
  Analogous to Proposition~9.10 in~\cite{mazur_rubin_ds}, the set $\mathcal{P}_0$ has positive density as explained as follows:
  By assumption the Galois group of $k(E[p])/k(\mu[p])$ is isomorphic to $\SL_2(\FF_p)$, which cannot have a subextension that is cyclic of order $p$ over $k(\mu[p])$.
  Hence the extension $L/k(\mu[p])$ obtained by adjoining all $p$-th roots of $\Sigma$-units intersects $k(E[p])$ exactly in $k(\mu[p])$.
  To ask $\mathfrak p$ to lie in $\mathcal{P}$ imposes that the Frobenius in $L/k$ is trivial and to lie in $\mathcal{P}_0$ is to ask that the Frobenius in $\SL_2(\FF_p)$ does not have a fix point.
  These conditions are independent resulting in $\mathcal{P}_0$ having positive density by Chebotarev's density theorem.

  For each of the infinitely many $\mathfrak{p}$ in $\mathcal{P}_0$, the extension $k(\mathfrak{p})/k$ is ramified only at $\mathfrak{p}$.
  By construction, we have that $p\nmid \mathfrak{p}$ and $E$ has good reduction at $\mathfrak{p}$.
  Therefore Proposition~\ref{easy_d_prop} applies.
  This implies that $D_{K/k}$ is trivial as $p\nmid \#\tilde E(\FF_\mathfrak{p})$ and $f_v=1$ for all places $v$ with $p\mid c_v$.
  Finally, Proposition~\ref{cyclic_items_prop}~\ref{pp2_item} shows that $\rk E(K) = 0$ and hence $M=0$ for all $K=k(\mathfrak{p})$ for $\mathfrak{p}\in \mathcal{P}_0$.
\end{proof}

\begin{prop}\label{cyclic_pos2_prop}
  Let $E/k$ be an elliptic curve and $p$ and odd prime.
  Suppose that there are more than $\rk E(k)$ primes $v$ such that $p\mid c_v$.
  Then there is a positive proportion of cyclic extensions $K/k$ of degree $p$ for which we have $\rk E(K)>\rk E(k)$ or $\Sha_K \neq 0$.
\end{prop}

\begin{proof}
  There is a positive proportion of cyclic extensions $K/k$ of degree $p$ such that all places $v$ with $p\mid c_v$ are inert.
  For such a $K$, the dimension of $D_{K/k}$ is larger than the rank of $E(k)$ by Corollary~\ref{cyclic_d_cor}.
  Then Proposition~\ref{cokeralpha_prop} implies that $\coker(\alpha)$ is non-trivial.
  This implies that $H^1(G,M) \neq 0$ or $\coker(\eta) \neq 0$ by Lemma~\ref{sha_descent_lem}.
  In the first case $M$ contains copies of $A$ and hence the rank of $E(K)$ is larger than the rank of $E(k)$.
  In the second case, $\Sha_K^G$ and hence $\Sha_K$ is non-trivial.
\end{proof}

Let us specialise to the case when $k=\QQ$.
Since the $L$-function $L(E,\chi,s)$ admits an analytic continuation for all $\chi$, it is easy to determine when the rank grows, that is when $\rk E(K)>\rk E(k)$.
Note that we can calculate the value $L(E,\chi,1)$ very quickly using modular symbols and if that value is non-zero, then the rank does not grow.
Under our assumption that $\Sha_K$ is finite, $L(E,\chi,1)=0$ implies that the rank grows.
Rank growth is relatively rare, especially for $p>5$, as expected by the conjectures made in~\cite{david_fearnley_kisilevsky_2, fearnley_kisilevsky_kuwata}.
Cases where the rank grows by more than $p-1$ are hard to find, but see Example~\ref{rk6_ex} below for such a case.

\begin{cor}
  If $p=3$ or $p=5$, then a positive proportion of $(E,K)$ where $E/\QQ$ is an elliptic curve, ordered by height, and $K$ is a cyclic extension of degree~$p$ and prime conductor, ordered by conductor, satisfy $E(K)=0$.
\end{cor}
\begin{proof}
  This is a consequence of Theorem~\ref{cyclic_pos_thm} and the results by Bhargawa and Shankar in~\cite{bhargava_shankar, bhargava_shankar_5} which show that a positive proportion of $E/\QQ$ have trivial $p$-Selmer group over $\QQ$ when $p=3$ or~$5$.
  The restriction on the Galois representation is negligible.
\end{proof}

For all practical purposes, we can consider the calculation of $\ord_{s=1} L(E,\chi,s)$, which leads to a proven upper bound for $r_K=\rk E(K)$, the calculation of $r_{\QQ}$ and $\Sha_\QQ$ as easy.
So is the determination of $D_{K/\QQ}$.
As a consequence, in most cases we can calculate $M$ effectively from our methods without having to do any point search or descent for $E$ over $K$.

\subsubsection{Examples}\label{ex_cyclic_subsec}

All elliptic curves in this list of examples are given by their Cremona label as in~\cite{cremona} and provided with a link to the \href{https://www.lmfdb.org}{lmfdb}~\cite{lmfdb}.
The computational results are obtained using SageMath~\cite{sage}.
The $p$-primary parts of Tate-Shafarevich groups over $\QQ$ are proven correct by the methods used in~\cite{shark}.

The conductor of a cyclic extensions $K/\QQ$ of degree $p$ is the smallest $m$ such that $K\subset \QQ(\zeta_m)$.
If $m$ is prime, then there is a unique such $K$ in $\QQ(\zeta_m)$, and hence we only need to specify $m$ to give $K$.

For a character $\chi$ of $K/\QQ$, seen as a Dirichlet character modulo $m$, we define the algebraic $L$-value
\[
  \mathcal{L}(E,\chi) = \sum_{a \bmod m}\bar\chi(a) \, \bigl[\tfrac{a}{m}\bigr]^{\chi(-1)}\ \in\QQ(\chi) = \QQ(\zeta_p)
\]
where $[\cdot]^\pm$ is the modular symbol attached to $E$, normalised as in~\cite{wiersema_wuthrich} and computed as in~\cite{numerical_wuthrich}.
Since $\mathcal{L}(E,\chi)$ is a non-zero multiple of $L(E,\chi,1)$, the vanishing of $\mathcal{L}(E,\chi)$ indicates that the rank of $E(K)$ is larger than the rank of $E(\QQ)$.
Conversely, if $\mathcal{L}(E,\chi)\neq 0$, then $E(K) = E(\QQ)$ under our assumption that $E(\QQ)[p]=0$.

\begin{enumerate}[label={\textbf{Example~\Alph*)}},ref=\Alph*,wide,itemindent=\labelsep+\labelwidth,parsep=0pt]
  \item Let $E$ be the elliptic curve with Cremona label \href{https://www.lmfdb.org/EllipticCurve/Q/67a1/}{67a1} and let $K$ be the \href{https://www.lmfdb.org/NumberField/7.7.594823321.1}{cyclic field of degree $p=7$ and conductor~$29$}.
  The curve has rank $0$ over $\QQ$ and $\coker(\alpha) = D_{K/\QQ}^{\vee}$ has dimension $1$ as the number of points in $\tilde E(\FF_{29})$ is divisible by $7$.
  Calculating $\mathcal{L}(E,\chi)\neq 0$ for a non-trivial character of $K$, proves that the rank over $K$ is still $0$.
  Therefore $M=0$ in this case.
  However, since $H^1(G,M)=0$, but $\coker(\alpha)\neq 0$, we have shown that $\Sha_K$ is non-trivial.
  In fact, the BSD conjecture over $K$ is equivalent to $\Sha(E/K)$ having $7^2\cdot 13^2$ elements.
  In our example, the space of $\Sha_K$ fixed by $G$ is $1$-dimensional.

  \item Similar to the previous example, we have a case with $M=\ZZ_p$, yet $\Sha_K\neq 0$.
  Take $E$ to be the curve \href{https://www.lmfdb.org/EllipticCurve/Q/37a1/}{37a1} of rank $1$ over $\QQ$ and $K$ to be the \href{https://www.lmfdb.org/NumberField/5.5.1982119441.1}{quintic field of conductor $211$} and $p=5$.
  Again $D_{K/\QQ}^{\vee}$ is of dimension~$1$ as $E(\FF_{211})$ is cyclic of order $5^2\cdot 3^2$.
  However, the generator $P$ of $E(\QQ)$ reduces to a point of order $45$ modulo $211$.
  This shows that $\coker(\alpha)\neq 0$.
  This together with $\mathcal{L}(E,\chi)\neq 0$ allows us to conclude that $\coker (\eta)=\Sha_K^G$ has dimension at least $1$.
  BSD says that $\Sha(E/K)$ is of order $5^4$.

  \item The curve \href{https://www.lmfdb.org/EllipticCurve/Q/681b3/}{681b3} has rank 0 over $\QQ$, but $\Sha_{\QQ}\neq 0$ for $p=3$.
  Consider the \href{http://www.lmfdb.org/NumberField/3.3.361.1}{cubic extension $K$ of conductor $19$}.
  Since $\mathcal{L}(E,\chi)=0$, the rank grows in this extension, which means that $M$ is a power of $A$.
  However $\coker(\alpha)$ is trivial.
  This implies that the capitulation kernel $C_{K/\QQ}$ is non-trivial.
  A $2$-descent reveals that $E(K)$ has rank $2$, which shows that $M\cong A$.
  Therefore only a $1$-dimensional subspace of $\Sha_{\QQ}$ capitulates in $K/\QQ$.
  Hence $\Sha_K$ is still non-trivial; it is of order $9$ according to BSD.

  \item Consider $p=5$, the curve \href{https://www.lmfdb.org/EllipticCurve/Q/21a1/}{21a1} over the \href{http://www.lmfdb.org/NumberField/5.5.2825761.1}{quintic extension $K$ of conductor 41}.
  The rank is $0$ over $\QQ$, but $\mathcal{L}(E,\chi)=0$ proves that the rank is positive over $K$.
  Since $\Sha_{\QQ}$ is trivial, but $\coker(\alpha)$ has dimension~$1$, we see that $M$ cannot contain more than one copy of $A$.
  Therefore $\rk E(K)=2$.
  Since the cokernel of $\eta$ must be trivial, we find that $\Sha_K^G=0$ and hence that $\Sha_K=0$ as $G$ is a $p$-group.
  Conjecturally $\Sha(E/K)$ is trivial.

  A similar argument works for $p=7$, the curve \href{https://www.lmfdb.org/EllipticCurve/Q/38b1/}{38b1} over the \href{http://www.lmfdb.org/NumberField/7.7.128100283921.1}{extension of conductor $71$}.

  \item The curve \href{https://www.lmfdb.org/EllipticCurve/Q/89a1/}{89a1} has rank $1$ over $\QQ$.
  As the algebraic $L$-value $\mathcal{L}(E,\chi)=0$ for the \href{http://www.lmfdb.org/NumberField/11.11.41426511213649.1}{degree~$11$ extension of conductor $23$}, i.e., $K=\QQ(\zeta_{23})^+$, the rank must grow.
  However $\alpha$ is surjective and $\Sha_{\QQ}=0$.
  Therefore $H^1(G,M)=0$, which implies that $M=\ZZ_p[G]$ is free.
  We can also conclude that $\Sha_K=0$.

  For rank~$1$ curves with rank growth, it is very frequent that $M$ is free.

   \item The curve \href{https://www.lmfdb.org/EllipticCurve/Q/130a3/}{130a3} has rank~$1$ over $\QQ$ and $\mathcal{L}(E,\chi)=0$ for the \href{http://www.lmfdb.org/NumberField/3.3.1849.1}{cubic field $K=\QQ(\alpha)$ of conductor~$43$} with $\alpha^3 + \alpha^2 - 14\,\alpha + 8 = 0$.
   The analytic rank tells us that $r_K=3$ and we know that $\Sha_{\QQ}$ is trivial.
   However, $\coker(\alpha) $ has $p=3$ elements.
   Hence $M$ contains at most one copy of $A$, but we cannot decide at this point whether $M$ is $\ZZ_p\oplus A$ or $\ZZ_p[G]$.
   In this case, we actually calculate the points in $E(K)$.
   One finds that there is a point $P\in E(K)$ with $x$-coordinate $\tfrac{1}{16}\,(-48 -132\,\alpha + 33\,\alpha^2 )$.
   The usual saturation shows that $\ZZ_p[G]\, P$ has index $p$ in $M$, which already tells us that $M\cong \ZZ_p\oplus A$.
   Alternatively, one can calculate the matrix of how a non-trivial element $\sigma\in G$ acts on the saturated group and calculate the cohomology group $H^1(G,M)$ directly.

 \item \label{rk6_ex}
  The most interesting example is the curve \href{https://www.lmfdb.org/EllipticCurve/Q/5692/a/1}{5692a1} with $K$ the \href{https://www.lmfdb.org/NumberField/3.3.81.1}{cubic extension of conductor 9}, which already appears in~\cite{wuthrich_fine}.
  The curve has rank $2$ over $\QQ$, generated by $P_1=(0,5)$ and $(2,-1)$.
  The group $D_2$ is cyclic of order $3$ and the image of the norm map identifies with the points of good reduction as the reduction type is IV.
  As the Tamagawa number at the only other bad prime is coprime to $3$, the only other non-trivial $D_v$ is for $v=3$.
  Here the reduction is good ordinary with $6$ elements in the reduction $\tilde E(\FF_3)$.
  We are in the situation of Proposition~\ref{p_ram_prop}.
  The map $\tilde E(\FF_3)[3] \to \hat E(3\ZZ_3)/\hat E(9\ZZ_3)$ can shown to be surjective, which implies that $D_3$ is cyclic of order $3$ as it identifies with $\tilde E(\FF_3)/3 \tilde E(\FF_3)$.

  Since $P_1 - P_2$ has good reduction at $2$ and the reduction at $3$ is of order $2$, the map $E(\QQ)\hatox\ZZ_3 \to D_{K/\QQ}$ is not surjective.
  We conclude that $\coker(\alpha)$ is of dimension~$1$.
  This only reveals that $M$ contains at most one copy of $A$.
  Since the rank of $M$ can be determined to be $6$, we can already conclude that $M$ must contain at least one factor of $\ZZ_3[G]$.

  To complete the calculation and prove that $M\cong \ZZ_3\oplus A\oplus \ZZ_3[G]$ seems to require once more the calculation of $M$ and the action of $G$ explicitly on it.
  With the explicit basis in~\cite{wuthrich_fine}, this is not difficult to do.

\end{enumerate}

\subsection{Dihedral group}
We suppose now that $G$ is the dihedral group $D_p$ of order~$2p$.
The $p$-Sylow subgroup is $N$ and pick one subgroup $H$ of order $2$.
There is a unique non-trivial Brauer relation $\Theta = \Theta_2=1 - 2\cdot H - N + 2 \cdot G$.
Write $F$ for the field fixed by $N$ and $L$ for the field fixed by our chosen $H$.
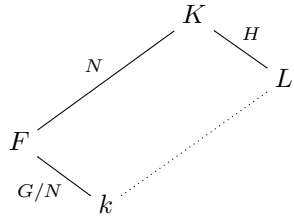
\[
\xymatrix@C-1ex@R-3ex{
  & &K\ar@{-}[dr]^H\ar@{-}[ddll]_N& \\
  & & &L\ar@{.}[ddll]\\
 F\ar@{-}[dr]_{G/N}& & &\\
  &k& & }
\]

Recall the classification of indecomposable $\ZZ_p[G]$-lattices: $\ZZ_p$, $\twist{\ZZ}_p := \ZZ_p\{1\}$, $A$, $\twist{A} := A\{1\}$, $B$ and $\twist{B}:=B\{1\}$.
From Theorem~\ref{torzewski_thm} we know that the following will determine $M=E(K)\hatox \ZZ_p$ completely:
\begin{itemize}[nosep]
  \item $M\otimes \QQ$ as a $\QQ[G]$-module;
  \item $H^1(N,M)$ as a $\FF_p[G/N]$-module;
  \item $s(M)= \ord_p(\mathcal{C}_{\Theta}(M))$.
\end{itemize}

 The last entry in the list above can be replaced by $\dim_{\FF_p} \iota(M)$.
 However, note that these invariants are not all easy to determine.
 The ranks could, at least conjecturally, be determined using the order of vanishing of twisted $L$-functions.
 The cohomological term $H^1(N,M)$ appears in the exact sequence~\eqref{shadescent_seq}.
 Finally, both $s(M)$ and $\iota(M)$ seem hard to evaluate without actually calculating $M$, except for the parity of $s(M)$ by Corollary~\ref{parity_cor}.

In the (frequent) case that the rank is small, we need less information to determine $M$.

\begin{prop}
  \begin{itemize}
    \item If $\rk E(F)=0$, then $M$ is determined by $H^1(N,M)$ as a $\FF_p[G/N]$-module.
    \item If $\rk E(F)=1$, then $M$ is determined by $H^1(N,M)$ as a $\FF_p[G/N]$-module, $\rk E(k)$ and the parity of $s(M)$.
  \end{itemize}
\end{prop}
\begin{proof}
  We use Table~\ref{dihedral_table}.
  If $r_F=0$, then $M$ is a direct sum of copies of $A$ and $\twist{A}$.
  Since they have distinct $H^1(N,M)$, that group is enough to determine $M$.

  If $r_F=1$, then $M$ is a direct sum of copies of $A$, $\twist{A}$ and one copy of either $\ZZ_p$, $\twist{\ZZ}_p$, $B$ or $\twist{B}$.
  Again $H^1(N,M)$ determines the number of $A$ and $\twist{A}$ that appear.
  If $r_k=1$, then there is an extra copy of either $\ZZ_p$ or $B$.
  Since $s(\ZZ_p)\equiv 1$ and $s(B)\equiv 0 \pmod{2}$, the parity of $s(M)$ suffices to determine $M$.
  If $r_k=0$, the same argument works with $\twist{\ZZ}_p$ and $\twist{B}$.
\end{proof}

 As $F/k$ is a quadratic extension, there is a quadratic twist $\twist{E}/k$ of $E$ associated to $F/k$.
 Since $p$ is odd, $\twist{E}$ also satisfies Assumption~\ref{tors_assu}.
 Many invariants that we might have to calculate over $F$ can be calculated over $k$ instead using $\twist{E}$.
 First of all $r_F = \rk E(F) = \rk E(k) + \rk\twist{E}(k)$.
 Since $p$ is odd, we also have $\Sha_F = \Sha_k \oplus \twist{\Sha}_k$, where $\twist{\Sha}_k = \Sha(\twist{E}/k)[p^\infty]$.
 Therefore, we also have $C_{K/F} = C_{K/k} \oplus \twist{C}_{K/k}$ with $\twist{C}_{K/k}$ the capitulation kernel for $\twist{E}$.

\begin{lem}
  $H^1(N,M)$ as a $G/N$-module is determined by the abelian groups $H^1(G,M)$ and $H^1\bigl(G,\twist{M}\bigr)$ where $\twist{M}$ is $\twist{E}(K)\otimes\ZZ_p = M\otimes \twist{\ZZ}_p$.
\end{lem}

\begin{proof}
  The $\FF_p$-vector space $H^1(N,M)$ splits into a $+1$ eigenspace and a $-1$ eigenspace with respect to the action by $G/N$.
  The $+1$ eigenspace is isomorphic to the $G/N$-invariants of $H^1(N,M)$, which is isomorphic to $H^1(G,M)$ by the restriction map.
  Twisting by $\twist{\ZZ}_p$, we obtain that the $-1$ eigenspace is $H^1(G,\twist{M})$.
\end{proof}

\begin{prop}\label{torzewski2_prop}
   $M=E(K)\hatox \ZZ_p$ is completely determined by
\begin{itemize}[nosep]
  \item $r_k=\rk E(k)$, $\twist{r}_{k} = \rk \twist{E}(k)$ and $\rk E(L)$;
  \item $H^1(G,M)$ and $H^1(G,\twist{M})$;
  \item $s(M)= \ord_p(\mathcal{C}_{\Theta}(M))$.
\end{itemize}
\end{prop}
\begin{proof}
  The structure of $M\otimes \QQ$ as a $\QQ_p[G]$ is determined by $r_k$, $\rk E(F) = r_k+\twist{r}_k$ and $\rk E(L)$.
  Together with the previous lemma, this theorem is now a reformulation of Theorem~\ref{torzewski_thm}.
\end{proof}

\begin{proof}[Proof of Theorem~\ref{ex1_thm}]
  The assumptions in Theorem~1 imply that $D_{K/k}$ is trivial using Proposition~\ref{easy_d_prop} as there are no places with $p\mid e_v$ and all places with $p\mid f_v$ have $p\nmid c_v$.
  Lemma~\ref{sha_descent_lem} together with the assumption that $\Sha_k$ is trivial, let us conclude that $H^1(G,M)$ is trivial.
  Also $H^1(G,\twist{M})$ is trivial for the same reason applied to $\twist{E}$.
  This implies that neither $A$ nor $\twist{A}$ can appear in $M$.
  Then $r_{k}+\twist{r}_k\leq 1$, implies that $M$ is isomorphic to a single copy of $\ZZ_p$, $\twist{\ZZ}_p$, $B$ or $\twist{B}$, unless $r_k=\twist{r}_k=0$, in which case $M=0$.
  If $r_k=1$, it is either $\ZZ_p$ or $B$, and if $\twist{r}_k=1$, it is  $\twist{\ZZ}_p$ or $\twist{B}$.
  The parity of $s(M)$ distinguishes the two possibilities in both cases, and that parity can be calculated using only local information for $E$ over $K$.
\end{proof}

\begin{lem}\label{sm_lem}
  Suppose that $p>3$ and that no place of additive reduction ramifies in $K/k$.
  Let $v_1$ be the number of places in $k$ such that $E$ has split multiplicative reduction and such that $K$ contains a single ramified place above $v$.
  Let $v_2$ be the number of places in $k$ such that $E$ has non-split multiplicative reduction and there is a unique place above $v$ with ramification index $p$.
  Then $s(M)\equiv v_1+v_2\pmod{2}$.
  In particular, $s(M)$ is even if there is no place of bad reduction that ramifies in $K/F$.
\end{lem}

\begin{proof}
  By Corollary~\ref{parity_cor}, we need to calculate the contribution at each bad place $v$ in $k$ to the $p$-adic valuation of $C(E/K)/C(E/F)$.

  For additive places the contribution is an even power of $p$:
  Since $p>3$, the Tamagawa number is not divisible by $p$, and since $K/k$ is unramified at this place the quantity $u_v$ does no change.
  Hence local term is $u_v^{p-1}$ or $1$ depending whether there are $p$ places above each place in $F$ or only one.

  For multiplicative places, this is calculated in the table in Section~3.1 in~\cite{bartel_psel}.
\end{proof}

An interesting application connects our investigation to the discussion of the ``minimalist conjecture'' in~\cite{kellock_dokchitser}.
Recall that we are still assuming that all Tate-Shafarevich groups have finite $p$-primary parts.
It is conjectured that a large proportion of elliptic curves have trivial $\Sha_{\QQ}$ by Delaunay~\cite{delaunay}.
For instance, this proportion is believed to be larger than $75.7$\%\ of all curves if $p\geq 3$.
We will make an assumption of this nature in the next theorem.
In the following, a ``large family'' of elliptic curves is understood in the sense of the definition at the start of Section~4 in~\cite{bhargava_shankar_5}, which means it is obtained by imposing relatively mild local conditions.

\begin{thm}\label{dihedral_pos_thm}
  Let $K/\QQ$ be a dihedral extension with $G=D_p$ for a prime $p>3$ such that $p\nmid e_p$.
  Suppose that a proportion of more than $\tfrac{17}{24}$ of elliptic curves $E/\QQ$ (when ordered by height) in any large family satisfy $\Sha(E/\QQ)[p]=0$.
  Then there is a positive proportion of elliptic curves $E/\QQ$, when ordered by height, such that $E(K)\hatox \ZZ_p$ is one of the following five $\ZZ_p[G]$-lattice:
  \[ 0,\qquad B,\qquad \twist{B},\qquad \ZZ_p\oplus \twist{\ZZ}_p,\qquad \text{and}\qquad B\oplus \twist{B} \cong \ZZ_p[G].\]
  The rank of $E(\QQ)$ and of the quadratic twist  $\twist{E}(\QQ)$ determine the case, except for the last two cases.
\end{thm}

We cannot conclude that the most frequent $\ZZ_p[G]$-module structures among the curves with non-trivial $\Sha_{\QQ}$ are the same as in the above theorem, but this could be true.
One would have to understand the frequency with which non-trivial elements in Tate-Shafarevich groups capitulate in $K$.

When tensoring the displayed formula in the theorem by $\CC_p$, one falls onto the ``minimalist conjecture'', except in the case that $M\cong B\oplus \twist{B} \cong\ZZ_p[G]$.
With our methods we cannot determine that the case $\ZZ_p\oplus\twist{\ZZ}_p$ is more frequent than $\ZZ_p[G]$.
Apart from that, the theorem is good evidence for the minimalist conjecture.

\begin{proof}
 We are going to use the work of Bhargava and Shankar~\cite{bhargava_shankar_5}, which already proves that a positive proportion of elliptic curve have rank $0$ or $1$.
 Consider the set of elliptic curves in
 \[
   \Bigl\{ E/\QQ \Bigm\vert (N,\Delta_K)=1,\ \ell^p\nmid \Delta_E\text{ for all primes $\ell$},\ p\nmid \#\tilde E(\FF_v)\text{ for all prime $v\in S_{\mathrm{r}}$} \Bigr\}
 \]
 where $S_{\mathrm{r}}$ is the set of all places $v$ such that $p\mid e_v$.
 This is a large family as defined in~\cite{bhargava_shankar_5}.

 Therefore more than $\tfrac{19}{24}$ of such elliptic curves have rank either $0$ or $1$ as shown in Proposition~38 in~\cite{bhargava_shankar_5}.
 Hence at least $1-2\cdot(1-\tfrac{19}{24})= \tfrac{7}{12}$ of all curves in the set have rank smaller than $2$ and their twist corresponding to the quadratic extension in $K/\QQ$ also have rank smaller than $2$.
 By assumption, strictly less than $2\cdot (1-\tfrac{17}{24}) =  \tfrac{7}{12}$ of  curves in the large set have $\Sha_{\QQ}$ or $\twist{\Sha}_{\QQ}$ non-trivial.
 This proves that a positive proportion of curves in the above set will now have rank either $0$ or $1$ for $E$ and its twist and trivial $\Sha_\QQ$ and $\twist{\Sha}_\QQ$.
 We may exclude the elliptic curves with a rational $p$-torsion point without harming this.

 Let $E$ be a such a curve.
 Since $\ell^p\nmid \Delta$, we see that $c_\ell$ cannot be divisible by $p$.
 Together with the condition that no bad prime ramifies, that $p\nmid e_p$ and the condition $p\nmid \# \tilde E(\FF_v)$ for all $v\in S_{\mathrm{r}}$, we deduce that $D_{K/\QQ}=0$ by Proposition~\ref{easy_d_prop}.
 Together with $\Sha_\QQ=\twist{\Sha}_\QQ=0$, we know now that $H^1(G,M)=0$ by Lemma~\ref{sha_descent_lem}.
 Therefore $M$ is a direct sum $\ZZ_p^a \oplus \twist{\ZZ}_p^b \oplus B^e \oplus \twist{B}^f$ with $r=\rk E(\QQ) = a+e\leq 1$ and $\twist{r}=\rk\twist{E}(\QQ) =b+f\leq 1$.

 By Lemma~\ref{sm_lem}, the quantity $a+b$ must be even as we are in the case that no prime of bad reduction ramified in $K/\QQ$.

 If $r=\twist{r}=0$, then $M$ must be $0$.
 If $r=1$, but $\twist{r}=0$, then $b=f=0$, which implies that $a=0$ since $a+b$ must be even, and hence $e=1$.
 Similar if $r=0$ and $\twist{r}=1$, we get $a=b=e=0$ and $f=1$.

 Finally, if $r=1$ and $\twist{r}=1$, then $a+e=1$ and $b+f=1$ and $a+b$ is even.
 This leads to two possibilities, namely $M\cong \ZZ_p\oplus \twist{\ZZ}_p$ or $M\cong B\oplus \twist{B} \cong \ZZ_p[G]$.
\end{proof}

 The four cases can also be determined by root numbers if one admits the parity conjectures.
 These root numbers can be calculated under our assumption (Theorem~2.15 in~\cite{kellock_dokchitser} as done in their Example~4.11):
 If $\Delta_F$ is the fundamental discriminant of $F$, then the root number of $E$ twisted by any of the irreducible $2$-dimensional $\CC[G]$-modules is $z = \sign(\Delta_F) \cdot \bigl( \Delta_F/N)$ where the second factor is the Jacobi symbol.
 The root number for $E$ twisted by the non-trivial quadratic character is $z$ times the root number of $E/\QQ$.
 The fact that the product of the three root numbers is $1$ is now equivalent to the result in Lemma~\ref{sm_lem} under the parity conjecture.

\subsection{Examples}\label{ex_dihedral_subsec}

For the following examples, we will always take the same $D_3$-extension.
Let \href{http://www.lmfdb.org/NumberField/3.1.140.1}{$L$ be the field generated by $\alpha$ with $\alpha^3 + 2\alpha -2=0$} and let \href{http://www.lmfdb.org/NumberField/6.0.686000.1}{$K$ be its Galois closure}.
The quadratic field inside $K$ is \href{http://www.lmfdb.org/NumberField/2.0.35.1}{$F=\QQ\bigl(\sqrt{-35}\bigr)$}.

There is a unique ramified prime above $2$ in $K$ and it has ramification index $3$.
Above $5$ and $7$ there are three primes with ramification index $2$.
The prime $3$ is unramified with residue degree $3$.

 Let \[M={\ZZ}_p^a \oplus \twist{\ZZ}_p^b \oplus A^c \oplus \twist{A}^d \oplus B^e \oplus \twist{B}^f\]
 and we try to determine the unknown $a,b,c,d,e,f$ from $r_\QQ=a+e$, $r_F= a+e+b+f$, $r_L = a+e+\tfrac{p-1}{2}(c+d+e+f)$, $-a+b-c+d\equiv \ord_p(C(E/K)/C(E/F))\pmod{2}$, $H^1(N,M) = \FF_p^c \oplus \twist{\FF}_p^d$.

\begin{enumerate}[resume*]
\item
  We take the \href{http://www.lmfdb.org/EllipticCurve/Q/82a1/}{curve 82a1} whose rank over $\QQ$ is $1$ and it is also $1$ over $F$ as the twist $\twist{E}$ has rank $0$.
  Therefore $b=f=0$.
  For all places $v\neq 2$, Proposition~\ref{notwild_prop} implies that $D_v=0$.
  Let $\mathfrak{p} = (2)$ be the prime above $2$ in $F$.
  We can determine $D_{\mathfrak p}$ for the extension $K/F$ using Proposition~\ref{non_p_totram_prop} as $K/F$ is totally ramified at $\mathfrak{p}$ and $E$ has split multiplicative reduction at this place.
  The quantity $u$ turns out to be odd and hence $D_{\mathfrak{p}}$ is cyclic of order $3$.
  However the rational point $P=(0,0)\in E(F)$ reduces to a non-singular point that is not in the formal group.
  Therefore $\coker(\alpha_{K/F})$ is trivial.
  It follows that $\alpha$ for $K/\QQ$ is also surjective.

  Since $\Sha_{\QQ}$ and $\twist{\Sha}_{\QQ}$ are trivial, we conclude that $c=d=0$.
  We are left with two possibilities, either $\ZZ_p$ or $B$.
  However Corollary~\ref{parity_cor} can be used now to show that $a$ is odd, since $C(E/K)/C(E/F)=3$.

  Therefore $M\cong\ZZ_p$ and we obtained this information with local information and information about $E$ and $\twist{E}$ over $\QQ$ only.
  For this particular curve it is not much effort to verify that $r_L=1$ with a $2$-descent, which confirms this result.

\item\label{I_ex}
  The next curve we take is \href{http://www.lmfdb.org/EllipticCurve/Q/14a3/}{14a3} which has rank $0$ over $\QQ$, but rank $1$ over $F$.
  Again $D_{K/\QQ}$ is reduced to $D_2$.
  Over $F$, the curve has split multiplicative reduction with Tamagawa number $18$.
  As in the above example $D_{\mathfrak{p}}$ is cyclic of order $3$, but this time the rational points map trivially to $D_{\mathfrak{p}}$.
  Therefore $\coker(\alpha_{K/F})$ is cyclic of order $3$.
  The same argument works for the twisted curve $\twist{E}$ over $\QQ$, showing that $\coker(\alpha_{K/F})$ is isomorphic to $\twist{\FF}_3$ as a $G/N$-module.
  Since the Tate-Shafarevich groups are trivial again, we know that $H^1(N,M)$ is either trivial or equal to~$\twist{\FF}_3$.

  This implies that $c=0$ and $d\leq 1$.
  The regulator constant yields $b\not\equiv d \pmod{2}$.
  We now have two possibilities left $d=1$ (and then $a=c=e=f=0$ and $f=1$) or $d=0$ (and then $b=1$ and $a=c=e=f=0$); so either $M\cong \twist{A}\oplus\twist{B}$ or $M\cong\twist{\ZZ}_p$.
  The fact that the $L$-function of $E$ twisted with the irreducible representation $\rho$ does not vanish at $s=1$ or, directly, a $2$-descent over $L$ confirms that $M\cong \twist{\ZZ}_p$.

  As a consequence, $\coker(\alpha_{K/F})$ having dimension $1$ now implies that $\Sha_K\neq 0$.
  We expect $\Sha(E/K)$ to have $9$ elements.

\item
  The \href{https://www.lmfdb.org/EllipticCurve/Q/322b1/}{curve 322b1} has rank $0$ over $F$ and hence $a=b=e=f=0$.
  The regulator constant tells us that $c\not\equiv d \pmod{2}$.
  Once again $\coker(\alpha_{K/F})=D_{K/F} =\twist{\FF}_p$ as a $G/N$-module very much like in the previous example as the reduction at $2$ is once more non-split multiplicative.
  Therefore $d\leq 1$ and $c=0$.
  We conclude that $M\cong \twist{A}$ without having to use any $L$-values or $2$-descents.

\item
  The situation is very similar for the \href{https://www.lmfdb.org/EllipticCurve/Q/158e1/}{curve 158e1} has also rank $0$ over $F$, but this time $\coker(\alpha_{K/F}) \cong {\FF}_p$ as a $G/N$-module since the reduction at~$2$ is split multiplicative.
  The argument as above will show that $M\cong {A}$.
  The difference between the two cases is that here $E(L)\oplus \tau E(L)$ will be equal to $E(K)$ while in the previous example it has index $p$.
  This can be checked by calculating the groups directly.

\item
  Finally, let us consider the \href{https://www.lmfdb.org/EllipticCurve/Q/37a1/}{curves 37a1} and the \href{https://www.lmfdb.org/EllipticCurve/Q/57a1/}{curve 57a1}.
  Both have rank $1$ over $\QQ$ and rank $2$ over $F$ and all Tate-Shafarevich groups in sight are trivial.
  For both curves the map $\alpha$ is surjective, which means that $H^1(N,M)$ is trivial, and all bad places are unramified, which implies that $s(M)$ is even.
  Therefore we are in the situation in Theorem~\ref{dihedral_pos_thm} where we had two options that we could not distinguish.
  However determining the group $E(L)$ in both cases, reveals that for the curve 37a1, we have $M=\ZZ_p\oplus \twist{\ZZ}_p$, while for the curve 57a1 it is $M=B\oplus \twist{B}$.
\end{enumerate}

Examples of $M\cong B$ or $M\cong \twist{B}$ can be found by Theorem~\ref{dihedral_pos_thm} or explicitly in~\cite{bmw2}.
More details for the above an verification using explicit points on $E(K)$ are done in~\cite{vavasour}.

\bibliographystyle{amsplain}
\bibliography{mwgm}

\end{document}